\documentclass{article}

\usepackage[preprint]{neurips_2021}

\usepackage[utf8]{inputenc} 
\usepackage[T1]{fontenc}    
\usepackage{hyperref}       
\usepackage{url}            
\usepackage{booktabs}       
\usepackage{amsfonts}       
\usepackage{nicefrac}       
\usepackage{microtype}      
\usepackage{xcolor}         
\usepackage{amsmath,amssymb,amsthm}
\usepackage{mathtools}
\usepackage{subfig}
\usepackage{float}
\usepackage{multirow}
\usepackage{graphicx}
\usepackage{algorithmic,algorithm}
\usepackage{tcolorbox} 
\usepackage{wrapfig}
\usepackage{caption}
\graphicspath{{figures}}

\usepackage{comment}

\DeclareMathOperator*{\argmin}{arg\,min}
\DeclareUnicodeCharacter{2082}{$_2$}

\newcommand{\miniboone}{\textsc{Miniboone}}

\newcommand{\bfv}{\mathbf{v}}

\newcommand{\bfx}{\boldsymbol{x}}
\newcommand{\bfy}{\boldsymbol{y}}

\newcommand{\bbR}{\mathbb{R}}
\def\tr{\operatorname{trace}}

\newtheorem{remark}{Remark}

\newcommand{\hf}{\frac{1}{2}}

\title{
Taming Hyperparameter Tuning in Continuous Normalizing Flows Using the JKO Scheme
}

\author{%
  Alexander Vidal\thanks{vidal@mines.edu} \\
  Dept. of Applied Mathematics and Statistics\\
  Colorado School of Mines \\
  \And
  Samy Wu Fung\\
  Dept. of Applied Mathematics and Statistics\\
  Dept. of Computer Science\\
  Colorado School of Mines\\
  \And 
  Luis Tenorio\\
  Dept. of Applied Mathematics and Statistics\\
  Colorado School of Mines\\
  \And
  Stanley Osher\\
  Dept. of Mathematics\\
  University of California, Los Angeles\\
  \And
  Levon Nurbekyan\\
  Dept. of Mathematics\\
  University of California, Los Angeles\\
}
\begin{document}

\maketitle

\begin{abstract}
 A normalizing flow (NF) is  a mapping that transforms  a chosen probability distribution to a normal distribution. Such flows are a common technique used for data generation and density estimation in machine learning and data science. The density estimate obtained with a NF requires a change of variables formula that involves the computation of the Jacobian determinant of the NF transformation. In order to tractably compute this determinant, continuous normalizing flows (CNF) estimate the mapping and its Jacobian determinant using a neural ODE. Optimal transport (OT) theory has been successfully used to assist in finding  CNFs by formulating them as OT problems with a soft penalty for enforcing the standard normal distribution as a target measure. A drawback of OT-based CNFs is the addition of a hyperparameter, $\alpha$, that controls the strength of the soft penalty and requires significant tuning.  We present JKO-Flow, an algorithm to solve OT-based CNF without the need of tuning $\alpha$. This is achieved by integrating the OT CNF framework into a Wasserstein gradient flow framework, also known as the JKO scheme. Instead of tuning $\alpha$, we repeatedly solve the optimization problem for a fixed $\alpha$ effectively performing a JKO update with a time-step $\alpha$. Hence we obtain a "divide and conquer" algorithm by repeatedly solving simpler problems instead of solving a potentially harder problem with large $\alpha$.
\end{abstract}

\section{Introduction}\label{sec: intro}
   A normalizing flow (NF) is a type of generative modeling technique that has shown great promise in applications arising in physics~\cite{noe2019boltzmann,brehmer2020madminer,carleo2019machine} as a general framework to construct probability densities for continuous random variables in high-dimensional spaces~\cite{papamakarios2019normalizing,kobyzev2019normalizing,rezende2015}.
    An NF provides a $\mathcal{C}^1$-diffeomorphism $f$ (i.e., a normalizing transformation) that transforms the density $\rho_0$ of an initial distribution $P_0$ to the density $\rho_1$ of the standard multivariate normal distribution $P_1$ -- hence the term "normalizing." Given such mapping $f$, the density $\rho_0$ can be recovered from the Gaussian density via the change of variables formula,
    \begin{equation}
        \log \rho_0(x) = \log \rho_1\left(f(x)\right) + \log|\det J_f(x)|,
        \label{eq: Jacobi_identity}
    \end{equation}
    where $J_f \in \bbR^{d \times d}$ is the Jacobian of $f$. Moreover, one can obtain samples with density $\rho_0$ by pushing forward Gaussian samples via $f^{-1}$.

    \begin{remark}
        Throughout the paper we slightly abuse the notation, using the same notation for probability distributions and their density functions. Additionally, given a probability distribution $P_0$ on $\mathbb{R}^d$ and a measurable mapping $f:\mathbb{R}^d \to \mathbb{R}^d$, we define the pushforward distribution of $P_0$ through $f$ as $(f\sharp P_0)(B)=P_0(f^{-1}(B))$ for all Borel measurable $B\subseteq \mathbb{R}^d$~\cite{villani2003topics}.
    \end{remark}
    
    There are two classes of normalizing flows: finite and continuous. A finite flow is defined as a composition of a finite number of $\mathcal{C}^1$-diffeomorphisms: $f = f_1 \circ f_2 \circ \cdots \circ f_n$. To make finite flows computationally tractable, each $f_i$ is chosen to have some regularity properties such as a Jacobian with a tractable determinant; for example, $J_{f_i}$ may have a triangular structure~\cite{papamakarios2017masked, baptista2021learning, zech2022sparse}.
     

 On the other hand, continuous normalizing flows (CNFs) estimate $f$ using a neural ODE of the form~\cite{chen2017continuous}:
    \begin{equation}
        \partial_t z(x,t) = v_\theta(z(x,t),t), \qquad z(x,0) = x, \qquad 0 \leq t \leq T,
        \label{eq: neural_ode}
    \end{equation}
    where $\theta$ are the parameters of the neural ODE.
    In this case, $f$ is defined as $f(x)= z(x,T)$ (for simplicity, we remove the dependence of $z$ on $\theta$).
    One of the main advantages of CNFs is that we can tractably estimate the log-determinant of the Jacobian using Jacobi's identity, which is commonly used in fluid mechanics (see, e.g., ~\cite[p.114]{villani2003topics}):
    \begin{equation}\label{eq:Jacobi}
        \begin{split}
        \partial_t \log |\det \nabla_x z(x,t)| &= \nabla_z \cdot v_\theta(z(x,t),t)
        = \tr\left( \nabla_z v_\theta (z(x,t),t)\right).
        \end{split}
    \end{equation}
 This is computationally appealing as one can replace the expensive determinant calculation by a more tractable trace computation of $\nabla_z v_\theta(z(x,t),t)$.
Importantly, no restrictions on $\nabla_z v_\theta(z(x,t),t)$ (e.g., diagonal or triangular structure) are needed; thus, these Jacobians are also referred to as ``free-form Jacobians''~\cite{grathwohl2019ffjord}.

The goal in training a CNF is to find parameters, $\theta$, such that $f=z(\,\cdot,T)$ leads to a good approximation of $\rho_1$ or, assuming $f$ is invertible, the pushforward of $\rho_1$ through $f^{-1}$ is a good approximation of $\rho_0$
~\cite{rezende2015,papamakarios2017masked,papamakarios2019normalizing,grathwohl2019ffjord}. Indeed, let $\widehat{\rho}_0$ be this pushforward density obtained with a CNF $f$; that is, $\widehat{\rho}_0=f^{-1}\sharp \rho_1$. We then minimize the Kullback-Leibler (KL) divergence from $\widehat{\rho}_0$ to $\rho_0$ given by 
\[
\min_{\theta}~ \mathbb{E}_{x \sim \rho_0} \log ( \rho_0(x)/\widehat{\rho}_0(x))= \min_{\theta}~ \mathbb{E}_{x \sim \rho_0}\left[\log \rho_0(x) - \log \rho_1(z(x,T)) - \ell(x,T)\right],
\]
where $\ell(x,T) = \log | \det \nabla z(x,T) |$. Dropping the $\theta$-independent term $\log \rho_0$ and using \eqref{eq: neural_ode} and \eqref{eq:Jacobi}, this previous optimization problem reduces to the minimization problem
\begin{equation}\label{eq: min_likelihood}
    \min_{\theta} \;\; \mathbb{E}_{x \sim \rho_0}\; C(x,T),
    \quad C(x,T):= -\log \rho_1(z(x,T))  - \ell(x,T)
\end{equation}
 subject to ODE constraints
	\begin{equation}
	    \partial_t 
	    \begin{bmatrix}
	        z(x,t)
	        \\
	        \ell(x,t)
	    \end{bmatrix}
	    = 
	    \begin{bmatrix}
        
	        v_\theta(z(x,t), t)
	        \\
	        \tr\left( \nabla_z v_\theta (z(x,t),t)\right)
	    \end{bmatrix}, 
	    \qquad 
	    \begin{bmatrix}
	        z(x,0)
	        \\
	        \ell(x,0)
	    \end{bmatrix}
	    =
	    \begin{bmatrix}
	     x
         \\
         0
	    \end{bmatrix}.
	    \label{eq: min_likelihood_constraints}
	\end{equation}

The ODE~\eqref{eq: min_likelihood_constraints} might be stiff for certain values of $\theta$, leading to extremely long computation times. Indeed, the dependence of $v$ on $\theta$ is highly nonlinear and might generate vector fields that lead to highly oscillatory trajectories with complex geometry. 

	Some recent work leverages optimal transport theory to find the CNF~\cite{finlay2020train,onken2021ot}. In particular, a kinetic energy regularization  term (among others) is added to the loss to ``encourage straight trajectories'' $z(x,t)$. That is, the flow is trained by minimizing the following instead of \eqref{eq: min_likelihood}:
	\begin{equation}
	    \begin{split}
        \min_{\theta} \;\; \mathbb{E}_{x \sim \rho_0} \; \int_0^T \frac{1}{2} \|v_\theta(z(x,t),t)\|^2 dt + \alpha C(x,T)
        \end{split}
        \label{eq: regularized_likelihood}
	\end{equation}
	subject to~\eqref{eq: min_likelihood_constraints}. The key insight in~\cite{finlay2020train,onken2021ot} is that~\eqref{eq: min_likelihood} is an example of a degenerate OT problem with a soft terminal penalty and without a transportation cost. The first objective term in~\eqref{eq: regularized_likelihood} given by the time integral is the transportation cost term, whereas $\alpha$ is a hyperparameter that balances the soft penalty and the transportation cost. Including this cost makes the problem well-posed by forcing the solution to be unique~\cite{villani2008optimal}. Additionally, it enforces straight trajectories so that~\eqref{eq: min_likelihood_constraints} is not stiff. Indeed~\cite{finlay2020train,onken2021ot} empirically demonstrate that including optimal transport theory leads to faster and more stable training of CNFs. Intuitively, we minimize the KL divergence \emph{and} the arclength of the trajectories. 
	
	Although including optimal transport theory into CNFs has been very successful~\cite{onken2021ot,finlay2020train,yang2019,zhang2018monge}, there are two key challenges that render them difficult to train. 
	First, estimating the log-determinant in~\eqref{eq: min_likelihood} via the trace in~\eqref{eq: min_likelihood_constraints} is still computationally taxing and commonly used methods rely on stochastic approximations~\cite{grathwohl2019ffjord,finlay2020train}, which adds extra error.
    Second, including the kinetic energy regularization requires tuning of the hyperparameter $\alpha$. Indeed, if $\alpha$ is chosen too small
    in \eqref{eq: regularized_likelihood}, then the kinetic regularization term dominates the training process, and the optimal solution consists of not moving, i.e., $f(x) = x$.
    On the other hand, if $\alpha$ is chosen too large, we return to the original setting where the problem is ill-posed, i.e., there are infinitely many solutions. Finally, finding an "optimal"  $\alpha$ is problem dependent and requires tuning on a cases-by-case basis.

\subsection{Our Contribution}
    We present JKO-Flow, an optimal transport-based algorithm for training CNFs without the need to tune the hyperparameter $\alpha$ in~\eqref{eq: regularized_likelihood}.
    Our approach also leverages fast numerical methods for exact trace estimation from the recently developed optimal transport flow (OT-Flow)~\cite{onken2021ot, ruthotto2020machine}.
    The key idea is to integrate the OT-Flow approach into a Wasserstein gradient flow framework, also known as the Jordan, Kinderlehrer, and Otto (JKO) scheme~\cite{jordan1998variational}. 
    Rather than tuning the hyperparameter $\alpha$ (commonly done using a grid search), the idea is to simply pick any $\alpha$ and solve a sequence of "easier" OT problems that gradually approach the target distribution. Each solve is precisely a gradient descent in the space of distributions, a Wasserstein gradient descent, and the scheme provably converges to the desired distribution for all $\alpha>0$~\cite{salim20Wasserstein}.  
    Our experiments show that our proposed approach is effective in generating higher quality samples (and density estimates) and also allows us to reduce the number of parameters required to estimate the desired flow. 

    Our strategy is reminiscent of debiasing techniques commonly used in inverse problems. Indeed, the transportation cost that serves as a regularizer in~\eqref{eq: regularized_likelihood} introduces a bias -- the smaller $\alpha$ the more bias is introduced (see, e.g., \cite{tenIP}), so good choices of $\alpha$ tend to be larger. 
    One way of removing the bias and the necessity of tuning the regularization strength is to perform a sequence of Bregman iterations~\cite{osher05iterative,burger05nonlinear} also known as nonlinear proximal steps. Hence our approach reduces to debiasing via Bregman or proximal steps in the Wasserstein space. In the context of CNF training, Bregman iterations are advantageous due to the flexibility of the choice for $\alpha$. Indeed, the resulting loss function is non-convex and its optimization tends to get harder for large $\alpha$. Thus, instead of solving one harder problem we solve several ``easier'' problems.


\section{Optimal Transport Background and Connections to CNFs}    

Denote by $\mathcal{P}_2(\mathbb{R}^d)$ the space of Borel probability measures on $\mathbb{R}^d$ with finite second-order moments, and let $\rho_0,\rho_1 \in \mathcal{P}_2(\mathbb{R}^d)$. The quadratic optimal transportation (OT) problem (which also defines the Wasserstain metric $W_2$) is then formulated as
\begin{equation}\label{eq:OT}
    W_2^2(\rho_0,\rho_1)=\inf_{\pi \in \Gamma(\rho_0,\rho_1)} \int_{\mathbb{R}^{2d}} \|x-y\|^2 d\pi(x,y),
\end{equation}
where $\Gamma(\rho_0,\rho_1)$ is the set of probability measures $\pi \in \mathcal{P}(\mathbb{R}^{2d})$ with fixed $x$ and $y$-marginals densities  $\rho_0$ and $\rho_1$, respectively. Hence the cost of transporting a unit mass from $x$ to $y$ is $\|x-y\|^2$, and one attempts to transport $\rho_0$ to $\rho_1$ as cheaply as possible. In~\eqref{eq:OT}, $\pi$ represents a \textit{transportation plan}, and $\pi(x,y)$ is the mass being transported from $x$ to $y$. One can prove that $(\mathcal{P}_2(\mathbb{R}^d),W_2)$ is a complete separable metric space~\cite{villani2003topics}. OT has recently become a very active research area in PDE, geometry, functional inequalities, economics, data science and elsewhere partly due to equipping the space of probability measures with a (Riemannian) metric~\cite{villani2003topics,villani2008optimal,peyre2018computational,santambrogio2015optimal}.

As observed in prior works, there are many similarities between OT and NFs~\cite{onken2021ot,finlay2020train,yang2020potential,zhang2018monge}. This connection becomes more transparent when considering the dynamic formulation of~\eqref{eq:OT}. More precisely, the Benamou-Brenier formulation of the OT problem is given by ~\cite{benamou2000computational}:
    \begin{equation}
        \begin{split} 
        \frac{T}{2} W_2^2(\rho_0, \rho_1) = \inf_{v, \rho} \;\; &\int_0^T \int_{\bbR^d} \frac{1}{2}\|v(x,t)\|_2^2 \rho(x,t) dx dt 
        \\
        \mbox{ s.t. } \;\; &\partial_t \rho(x,t) + \nabla \cdot (\rho(x,t) v(x,t)) = 0 
        \\
        & \rho(x,0) = \rho_0(x), \;\; \rho(x,T) = \rho_1(x).
        \end{split}
        \label{eq: benamou_brenier_formulation}
    \end{equation}

    Hence, the OT problem can be formulated as a problem of flowing $\rho_0$ to $\rho_1$ with a velocity field $v$ that achieves minimal kinetic energy. The optimal velocity field $v$ has several appealing properties. First, particles induced by the optimal flow $v$ travel in straight lines. Second, particles travel with constant speed. Moreover, under suitable conditions on $\rho_0$ and $\rho_1$, the optimal velocity field is unique~\cite{villani2003topics}.

    Given a velocity field $v$, denote by $z(x,t)$ the solution of the ODE
    \begin{equation*}
        \partial_t z(x,t) = v(z(x,t),t), \qquad z(x,0) = x, \qquad 0 \leq t \leq T.
    \end{equation*}
    Then, under suitable regularity conditions, we have that the solution of the continuity equation is given by $\rho(\cdot,t)=z(\cdot,t)\sharp \rho_0$.
    Thus the optimization problem in~\eqref{eq: benamou_brenier_formulation} can be written as
    \begin{equation}\label{eq:BB_Lagrangian}
        \begin{split} 
        \inf_{v} \;\; &\int_0^T \int_{\bbR^d} \frac{1}{2}\|v(z(x,t),t)\|_2^2 \rho_0(x) dx dt 
        \\
        \mbox{ s.t. } \;\; & \partial_t z(x,t) = v(z(x,t),t),~ z(x,0) = x,~ 
         z(\cdot,T)\sharp \rho_0=\rho_1.
        \end{split}
    \end{equation}

This previous problem is very similar to~\eqref{eq: min_likelihood} with the following differences:
\begin{itemize}
    \item the objective function in \eqref{eq: min_likelihood} does not have the kinetic energy of trajectories,
    \item the terminal constraint is imposed as a soft constraint in~\eqref{eq: min_likelihood} and as a hard constraint in~\eqref{eq:BB_Lagrangian}, and
    \item $v$ in~\eqref{eq: min_likelihood} is $\theta$-dependent, whereas the formulation in~\eqref{eq:BB_Lagrangian} is in the non-parametric regime.
\end{itemize}

So the NF~\eqref{eq: min_likelihood} can be thought of as an approximation to a degenerate transportation problem that lacks transportation cost. Based on this insight one can regularize~\eqref{eq: min_likelihood} by adding the transportation cost and arrive at~\eqref{eq: regularized_likelihood} or some closely related version of it~\cite{onken2021ot,finlay2020train,yang2020potential,zhang2018monge}. It has been observed that the transportation cost (kinetic energy) regularization significantly improves the training of NFs.

\section{JKO-Flow: Wasserstein Gradient Flows for CNFs}
    While the OT-based formulation of CNFs in~\eqref{eq: regularized_likelihood} has been found successful in some applications~\cite{onken2021ot,finlay2020train,yang2020potential,zhang2018monge}, a key difficulty arises in choosing how to balance the kinetic energy term and the KL-divergence, i.e., on selecting $\alpha$. This difficulty is typical in problems where the constraints are imposed in a soft fashion.
    Standard training of CNFs typically involves tuning for a ``large but hopefully stable enough'' step size $\alpha$ so that the KL divergence  term is sufficiently small after training.
    To this end, we propose an approach that avoids the need to tune  $\alpha$ by using the fact that the solution to~\eqref{eq: regularized_likelihood} is an approximation to a backward Euler (or proximal point) algorithm when discretizing the Wasserstein gradient flow using the Jordan-Kinderlehrer-Otto (JKO) scheme~\cite{jordan1998variational}. 
    The seminal work in \cite{jordan1998variational} provides a gradient flow structure of the Fokker-Planck equation using an implicit time discretization.
    That is, given $\alpha > 0$, density at $k^{\text{th}}$ iteration, $\rho^{(k)}$, and terminal density $\rho_1$,
    ~one finds
    \begin{equation}
    \begin{split}
        \rho^{(k+1)} =& \argmin_{\rho \in \mathcal{P}_2(\bbR^d)} \; \frac{1}{2\alpha} W_2^2(\rho, \rho^{(k)}) + KL(\rho||\rho_1)\\
        =&\argmin_{v} \; \frac{1}{\alpha}\int_0^1 \int_{\bbR^d} \frac{1}{2}\|v(z(x,t),t)\|_2^2 \rho_0(x) dx dt  + KL(z(\cdot,1)\sharp \rho^{(k)}||\rho_1)\\
        \mbox{ s.t. } \;\; & \partial_t z(x,t) = v(z(x,t),t),~ z(x,0) = x,~ z(\cdot,T)\sharp \rho_0=\rho_1,
    \end{split}
        \label{eq: JKO_update}
    \end{equation}
    for $k=0,1,\ldots$, and $\rho^{(0)} = \rho_0$.
    Here, $\alpha$ takes the role of a step size when applying a proximal point method to the KL divergence using the Wasserstein-2 metric, and $\{\rho^{(k)}\}$ provably converges to $\rho_1$~\cite{jordan1998variational,salim20Wasserstein}. Hence, repeatedly solving~\eqref{eq:BB_Lagrangian} with the KL penalty acting as a soft constraint yields an arbitrarily accurate approximation of $\rho_1$. 
    In the parametric regime each iteration takes the form
    \begin{equation}
        \argmin_{\theta} \;\; \mathbb{E}_{x \sim \rho^{(k)}} \; \int_0^T \frac{1}{2} \|v_\theta(x,t)\|^2 dt + \alpha C(x,T)
        \quad \text{subject to} \quad \eqref{eq: min_likelihood_constraints}.
        \label{eq: OT-Flow-based-JKO}
	\end{equation}
 Thus we solve a sequence of problems~\eqref{eq: regularized_likelihood}, where the initial density of the current subproblem is given by the pushforward of the density generated in the previous subproblem.
\begin{algorithm}[t]
    \caption{Proposed Algorithm}
    \label{alg: JKO-CNF}
    \begin{algorithmic}[1]
        \STATE{{\textbf{Input:}} Samples from $\rho_0$, step size $\alpha >0$, number of steps $K$}
        \STATE{Initialize $\theta_1$ at random}
        \FOR{$k=1,\ldots,K$}
            \STATE Solve for $\theta_{k}$ using samples by solving~\eqref{eq: OT-Flow-based-JKO}\
            \STATE Update distribution of samples using $v_{\theta_{k}}$ \\
        \ENDFOR
        \STATE{{\textbf{Output:}} saved weights $\theta_1,\ldots, \theta_K$}
    \end{algorithmic}
    \end{algorithm}
Importantly, our proposed approach \emph{does not require tuning $\alpha$}. Instead, we solve a sequence of subproblems that is guaranteed to converge to $\rho_1$~\cite{jordan1998variational} prior to the neural network parameterization; see Alg.~\ref{alg: JKO-CNF}.
	While our proposed methodology can be used in tandem with any algorithm used to solve~\eqref{eq: OT-Flow-based-JKO}, an important numerical aspect in our approach is to leverage fast computational methods that use \emph{exact} trace estimation in~\eqref{eq: min_likelihood_constraints}; this approach is called OT-Flow~\cite{onken2021ot}. 
    Consequently, we avoid the use of stochastic approximation methods for the trace, e.g., Hutchinson's estimator~\cite{avron2011,hutchinson1990stochastic,tenIP}, as is typically done in CNF methods~\cite{grathwohl2019ffjord,finlay2020train}.
	A surprising result of our proposed method is that it empirically shows improved performance even with fewer number of parameters (see Fig.~\ref{fig:checkerboardWidth}).
	
\section{Related Works}
    \noindent \textbf{Density Estimation.}
    One of the main advantages of NFs over other generative models is that they provide density estimates of probability distributions using~\eqref{eq: Jacobi_identity}. That is, we do not need to apply a separate density estimation technique after generating samples from a distribution, e.g., as in GANs~\cite{goodfellow2020generative}.
    Multivariate density estimation is a fundamental problem in statistics~\cite{silverman1986density,scott2015multivariate}, High Energy Physics (HEP)~\cite{cranmer2001kernel} and in other fields of science dealing with multivariate data. For instance,  particle physicists in HEP study possible distributions from a set of high energy data. Another application of density estimation is in confidence level calculations of particles in Higgs searches at Large Electron Positron Colliders (LEP)~\cite{opal2000search} and  discriminant methods used in the search for new particles~\cite{cranmer2001kernel}. 

    \noindent \textbf{Finite Flows.} Finite normalizing flows~\cite{tabak2013family,rezende2015,papamakarios2019normalizing,kobyzev2019normalizing} use a composition of discrete transformations, where specific architectures are chosen to allow for efficient inverse and Jacobian determinant computations. 
	NICE~\cite{dinh2014nice}, RealNVP~\cite{dinh2016density}, IAF~\cite{kingma2016improved}, and MAF~\cite{papamakarios2017masked} use either autoregressive or coupling flows where the Jacobian is triangular, so the Jacobian determinant can be tractably computed. GLOW~\cite{kingma2018glow} expands upon RealNVP by introducing an additional invertible convolution step. 
	These flows are based on either coupling layers or autoregressive transformations, whose tractable invertibility allows for density evaluation and generative sampling.
	Neural Spline Flows~\cite{durkan2019neural} use splines instead of the coupling layers used in GLOW and RealNVP. Using monotonic neural networks, NAF~\cite{huang2018neural} require positivity of the weights, which UMNN~\cite{wehenkel2019unconstrained} circumvent this requirement by parameterizing the Jacobian and then integrating numerically.

    \noindent \textbf{Continuous and Optimal Transport-based Flows.} 
    Modeling flows with differential equations is a natural and commonly used method~\cite{suykens1998,welling2011bayesian,neal2011mcmc,salimans2015markov,ruthotto2021introduction,huang2020convex}.
	In particular, CNFs model their flow via a neural ordinary differential equation  ~\cite{chen2017continuous,chen2018neural,grathwohl2019ffjord}.  Among the most well-known CNFs are FFJORD~\cite{grathwohl2019ffjord}, which estimates the determinant of the Jacobian by accumulating its trace along the trajectories, and the trace is estimated using Hutchinson's estimator~\cite{avron2011,hutchinson1990stochastic,tenIP}.  
 	To promote straight trajectories, RNODE~\cite{finlay2020train} regularizes FFJORD with a transport cost $L(\bfx,T)$. 
 	RNODE also includes the Frobenius norm of the Jacobian $\| \nabla \bfv \|_F^2$ to stabilize training. 
	The trace and the Frobenius norm are estimated using a stochastic estimator and report speedup by a factor of 2.8. 

	Monge-Ampère Flows~\cite{zhang2018monge} and Potential Flow Generators~\cite{yang2019} similarly draw from OT theory but parameterize a potential function instead of the dynamics directly.
	OT is also used in other generative models~\cite{sanjabi2018,salimans2018improving,lei2018geometric,lin2019fluid,avraham2019parallel,tanaka2019discriminator}.
	OT-Flow~\cite{onken2021ot} is based on a discretize-then-optimize approach~\cite{onken2020do} that also parameterizes the potential function. To evaluate the KL divergence, OT-Flow estimates the density using an \emph{exact} trace computation following the work of~\cite{ruthotto2020machine}. 
	
    \noindent \textbf{Wasserstein Gradient Flows.}
    Our proposed method is most closely related to~\cite{fan2021variational}, which also employs a JKO-based scheme to perform generative modeling. But a key difference is that~\cite{fan2021variational} reformulates the KL-divergence as an optimization over difference of expectations (See~\cite[Prop. 3.1]{fan2021variational}); this makes their approach akin to GANs, where the density cannot be obtained  without using a separate density estimation technique. 
	Our proposed method is also closely related to methods that use input-convex CNNs~\cite{mokrov2021large, bunne2022proximal, alvarez2021optimizing}.~\cite{mokrov2021large} focuses on the special case with KL divergence as objective function.~\cite{alvarez2021optimizing} solve a sequence of subproblems different from the fluid flow formulation presented in~\eqref{eq: OT-Flow-based-JKO}. They also require an end-to-end training scheme that backpropagates to the initial distribution; this can become a computational burden when the number of time discretizations is large.~\cite{bunne2022proximal} utilizes a JKO-based scheme to approximate a population dynamics given an observed trajectory and focus on applications in computational biology. 
	Other related works include natural gradient methods~\cite{nurbekyan2022efficient} and implicit schemes based on the Wasserstein-1 distance~\cite{heaton2020wasserstein}.
	
	

\section{Numerical Experiments} \label{sec: numerical_experiments}
We demonstrate the effectiveness of our proposed JKO-Flow on a series of synthetic and real-world datasets. As previously mentioned, we compute each update in~\eqref{eq: JKO_update} by solving~\eqref{eq: regularized_likelihood} using the OT-Flow solver~\cite{onken2021ot}, which leverages fast and exact trace computations. We also use the same architecture provided in~\cite{onken2021ot}.
Henceforth, we shall also call the traditional CNF approach the ``single-shot'' approach.

\newcommand{\rottextTwo}[1]{\rotatebox{90}{\parbox{22mm}{\small \centering#1}}}

\begin{figure}[t] 
    \centering
	\begin{tabular}{cccc}
		$\alpha = 1$ & $\alpha = 5$  & $\alpha = 10$  & $\alpha = 50$ 
        \\
        \multicolumn{4}{c}{\textbf{MMD}$^{\boldsymbol{2}}$ \textbf{Values, Single Shot}}
        \\
        3.58e-2 & 3.56e-3 & 1.42e-3 & 1.26e-3 \\
        \includegraphics[width=0.22\textwidth]{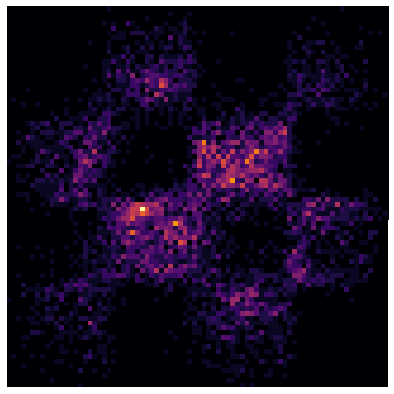}&
        \includegraphics[width=0.22\textwidth]{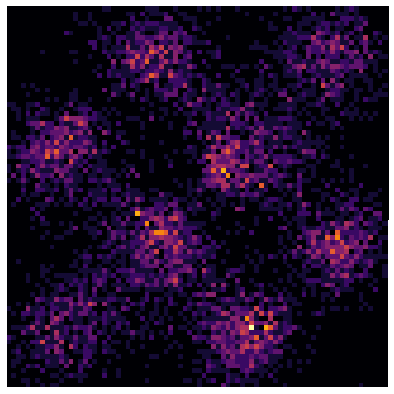}&
        \includegraphics[width=0.22\textwidth]{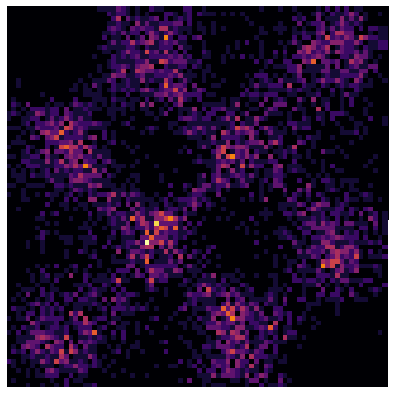}&
        \includegraphics[width=0.22\textwidth]{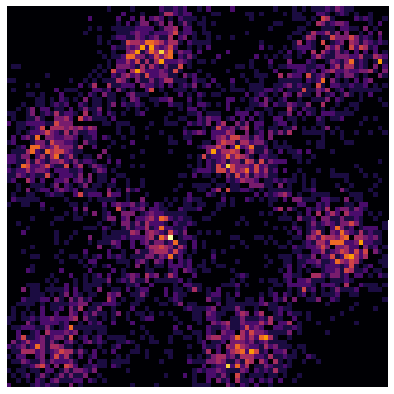}\\
               \\
        \multicolumn{4}{c}{\textbf{MMD$\boldsymbol{^2}$ Values, JKO-Flow (five iterations)}}
        \\
        4.90e-4 & 5.70e-4 & 6.40e-4 & 9.00e-4 \\
        \includegraphics[width=0.22\textwidth]{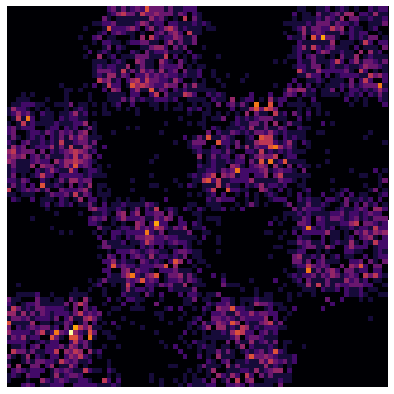}&
        \includegraphics[width=0.22\textwidth]{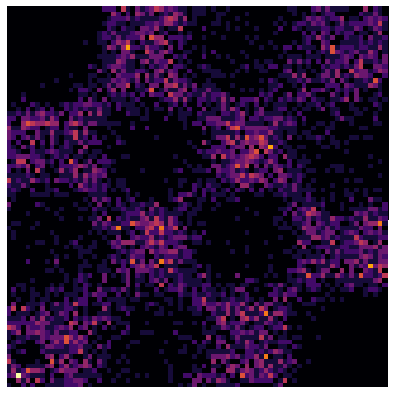}&
        \includegraphics[width=0.22\textwidth]{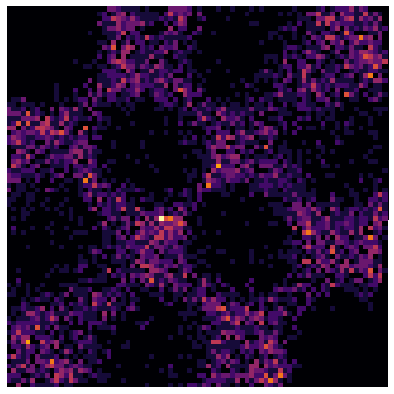}&
        \includegraphics[width=0.22\textwidth]{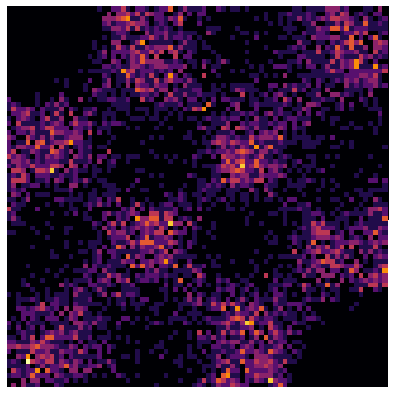} \\
	\end{tabular} 
	\centering 
	\\[2mm]
	\caption{{Checkerboard dataset:} Generated samples of $\widehat{\rho}_0$ using the standard one-shot approach (top row). Generated using our proposed JKO-Flow using five iterations (bottom row). Here, we use $\alpha  = $ 1, 5, 10, 50. JKO-Flow returns consistent results \emph{regardless of the value of $\alpha$}.}
  \label{fig:checkerboardalpha}
\end{figure} 

\begin{figure*}[t]
  \centering
  \setlength{\tabcolsep}{1pt}
  \renewcommand{\arraystretch}{0.5}
  \begin{tabular}{cccccccc}
  	\rottextTwo{Data\\$x$}
  	&
  	\includegraphics[width=0.13\textwidth]{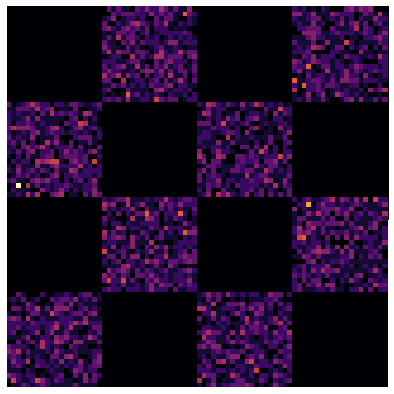}
  	&
  	\includegraphics[width=0.13\textwidth]{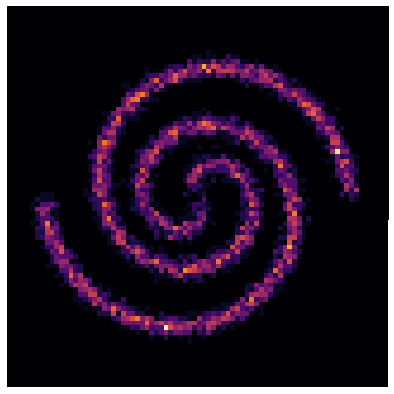}
  	&
  	\includegraphics[width=0.13\textwidth]{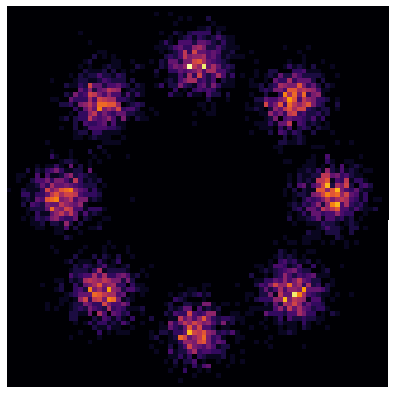}
  	&
  	\includegraphics[width=0.13\textwidth]{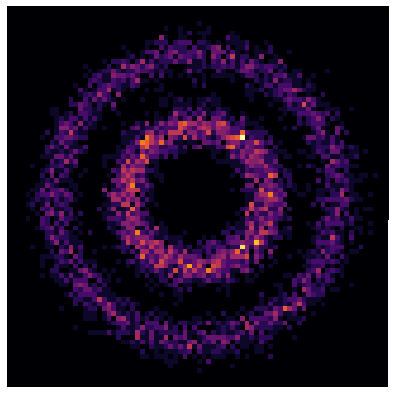}
  	&
  	\includegraphics[width=0.13\textwidth]{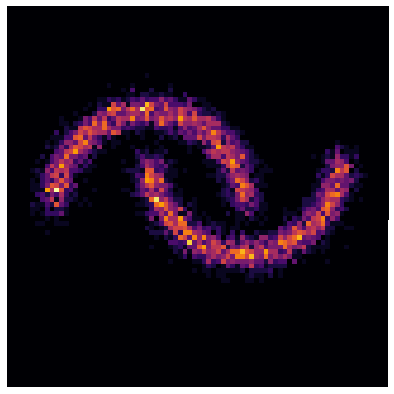}
  	&
  	\includegraphics[width=0.13\textwidth]{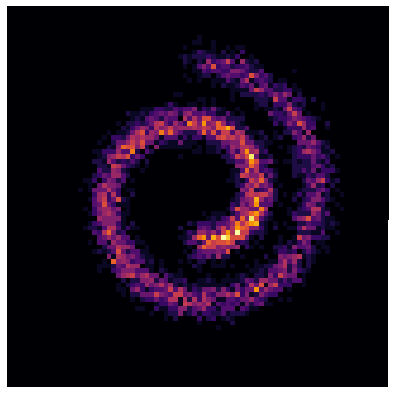}
  	&
  	\includegraphics[width=0.13\textwidth]{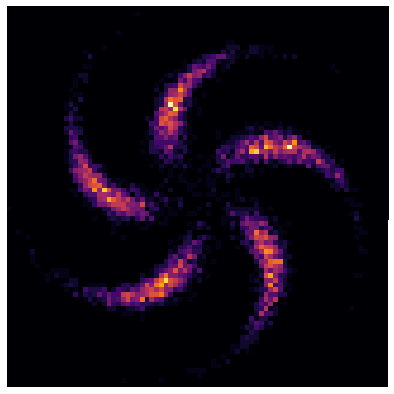}
  	\\
    \rottextTwo{Estimate\\$\rho_0$}
  	&
  	\includegraphics[width=0.13\textwidth]{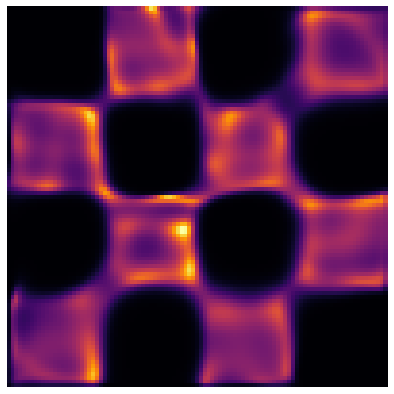}
  	&
  	\includegraphics[width=0.13\textwidth]{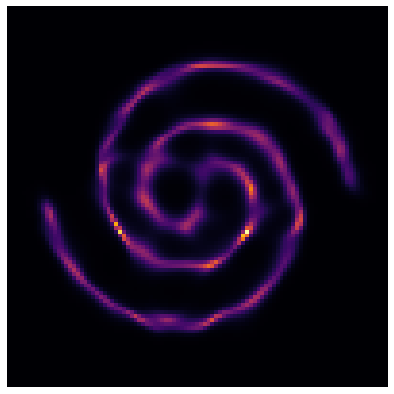}
  	&
  	\includegraphics[width=0.13\textwidth]{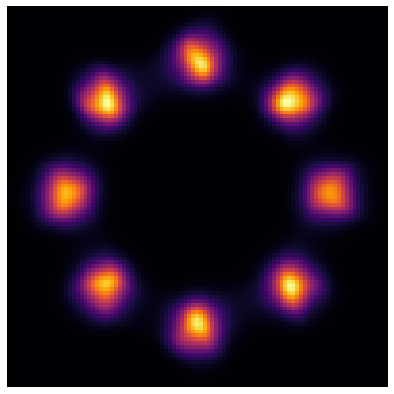}
  	&
  	\includegraphics[width=0.13\textwidth]{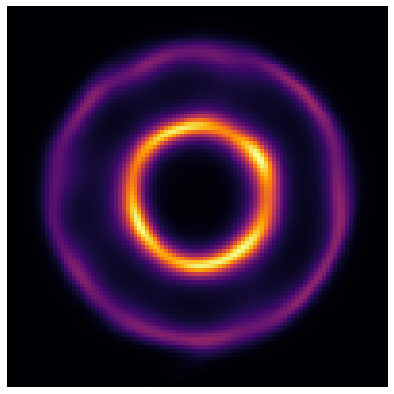}
  	&
  	\includegraphics[width=0.13\textwidth]{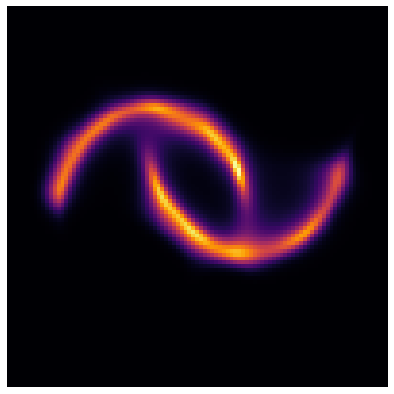}
  	&
  	\includegraphics[width=0.13\textwidth]{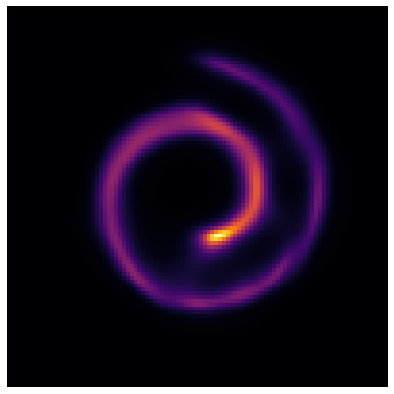}
  	&
  	\includegraphics[width=0.13\textwidth]{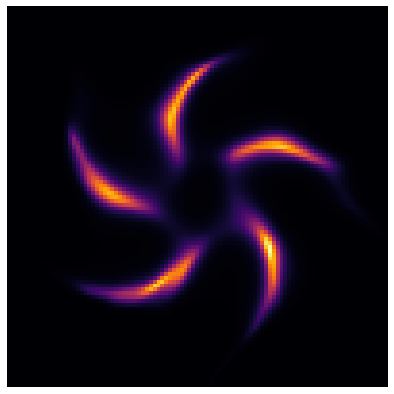}
  	\\
  	\rottextTwo{Generation\\$f^{-1}(y)$}
  	&
  	\includegraphics[width=0.13\textwidth]{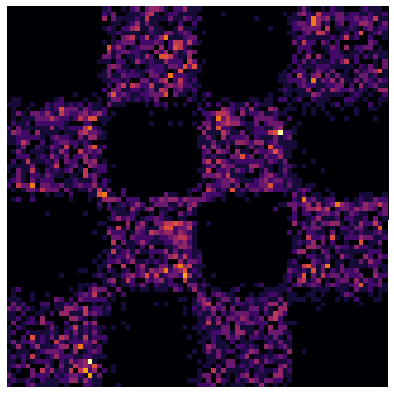}
  	&
  	\includegraphics[width=0.13\textwidth]{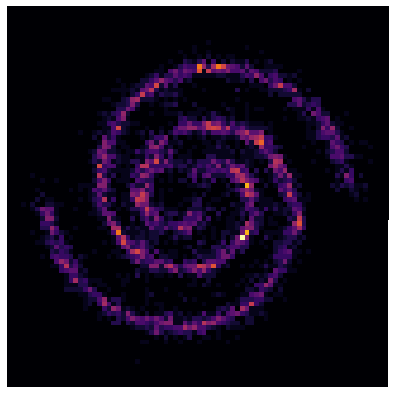}
  	&
  	\includegraphics[width=0.13\textwidth]{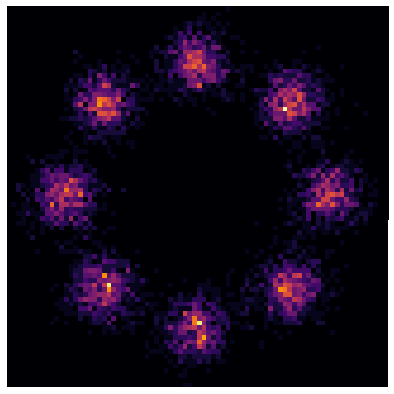}
  	&
  	\includegraphics[width=0.13\textwidth]{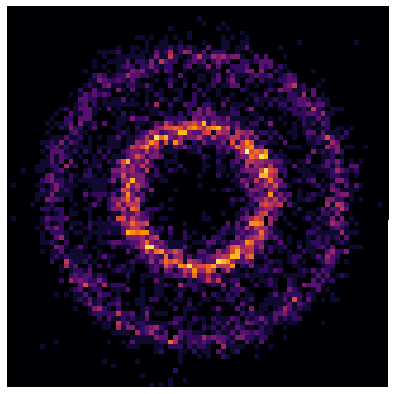}
  	&
  	\includegraphics[width=0.13\textwidth]{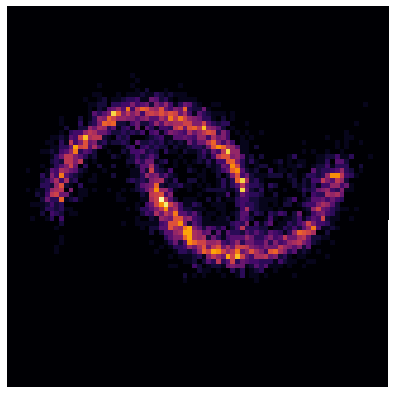}
  	&
  	\includegraphics[width=0.13\textwidth]{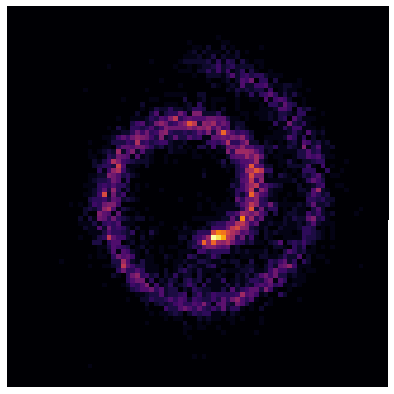}
  	&
  	\includegraphics[width=0.13\textwidth]{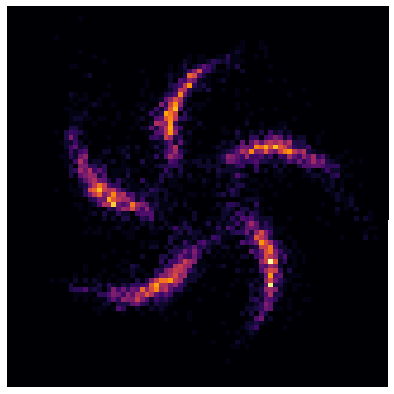}
  \end{tabular}
  \caption{Density estimation on 2D toy problems using five JKO-Flow iterations. \textbf{Top:} samples from the unknown distribution $\rho_0$. \textbf{Middle:} density estimate for $\rho_0$ computed by inverting the flow through the five iterations of JKO-Flow from $\rho_1$ via~\eqref{eq: neural_ode}.
  \textbf{Bottom:} samples generated by inverse JKO-Flow through five iterations where $y$ has density  $\rho_1$.}
  \label{fig:toyModels}
\end{figure*}

\noindent \textbf{Maximum Mean Discrepancy Metric (MMD).}  
Our density estimation problem requires approximating a density  $\rho_0$ by finding a transformation $f$
such that $f^{-1}\sharp \rho_1$ has density $\widehat{\rho}_0$ close to $\rho_0$ where $\rho_1$ is the standard Gaussian. However,
$\rho_0$ is not known in real-world density estimation scenarios, such as in physics applications, all we have are samples $X = \{x_i\}_{i=1}^n$
from the unknown distribution. 
Consequently, we use the observed samples $X$ and samples 
$\widehat{X}=\{\widehat{x}_j\}_{j=1}^m$, $\widehat{x}_j=f^{-1}(q_j)$,
generated by the CNF and samples $Q=\{q_j\}_{j=1}^m$ from $\rho_1$ to determine if their
corresponding distributions are close in some sense. To measure the discrepancy we use 
a particular integral probability metric \cite{zolotarev1976metric,rachev2013methods,muller1997integral} known as maximum mean discrepancy (MMD) defined as follows \cite{gretton2012mmd}: Let $x$ and $y$ be random vectors in $\mathbb{R}^d$ with distributions $\mu_x$ and $\mu_y$, respectively, and let $\mathcal{H}$ be a reproducing kernel Hilbert space of functions on $\mathbb{R}^d$ with Gaussian kernel (see~\cite{paulsen2016introduction} for an introduction)
\begin{equation}
    k(x_i, x_j) = \exp{\left(-\hf \|x_i - x_j\|^2\right)}.
\end{equation}
Then the MMD of $\mu_x$ and $\mu_y$ is given by
\[
\mathrm{MMD}_{\mathcal{H}}(\mu_x,\mu_y) = \sup_{\|f\|_\mathcal{H}\leq 1}\,|\,\mathbb{E}\,f(x) - \mathbb{E}\,f(y)\,|.
\]
It can be shown that $\mathrm{MMD}_\mathcal{H}$ defines a metric on the class of probability measures on $\mathbb{R}^d$ \cite{gretton2012mmd,fukumizu2007kernel}. 
The squared-MMD can be written in terms of the kernel as follows:
\[
\mathrm{MMD}^2_{\mathcal{H}}(\mu_x,\mu_y) =
\mathbb{E}\,k(x,x^\prime) + \mathbb{E}\,k(y,y^\prime)
-2\,\mathbb{E}\,k(x,y),
\]
where $x,x^\prime$ are iid $\mu_x$ independent of $y,y^\prime$ which are iid $\mu_y$.
An unbiased estimate of the squared-MMD based on the samples $X$ and $\widehat{X}$  defined above is given by \cite{gretton2012mmd}:
\[
\mathrm{MMD}^2_\mathcal{H}(X,\widehat{X}) =
\frac{1}{n(n-1)}\sum_{i\neq j} k(x_i,x_j)
+ \frac{1}{m(m-1)}\sum_{k\neq \ell} k(\widehat{x}_k,\widehat{x}_\ell)
-\frac{2}{nm}\sum_{i,\ell} k(x_i,\widehat{x}_\ell).
\]
Note that the MMD is not used for algorithmic training of the CNF, it is only used to compare the densities $\rho_0$ and $\widehat{\rho}_0$ based on the samples $X$ and $\widehat{X}$.

\noindent \textbf{Synthetic 2D Data Set.}  We begin by testing our method on  seven two-dimensional (2D) benchmark datasets for density estimation algorithms commonly used in machine learning \cite{grathwohl2019ffjord,wehenkel2019unconstrained}; see Fig.~\ref{fig:toyModels}.
We generate results with JKO-Flow for different values of $\alpha$ and for different number of iterations. 
We use $\alpha = 1, \, 5, \, 10, \,$ and $50$, and for each $\alpha$ we use the single shot approach $k=1$ and JKO-Flow with $k=5$ iterations from~\eqref{eq: JKO_update}.
Note that in CNFs, we are interested in estimating the density (and generating samples) from $\rho_0$; consequently, once we have the optimal weights $\theta^{(1)}, \theta^{(2)}, \ldots, \theta^{(5)}$, we must ``flow backwards'' starting with samples from the normal distribution $\rho_1$.
Fig.~\ref{fig:checkerboardalpha} shows that JKO-Flow outperforms the single shot approach for different values of $\alpha$. In particular, the performance for the single shot approach varies drastically for different values of $\alpha$, with $\alpha=1$ being an order of magnitude higher in MMD than $\alpha=5$.
On the other hand, JKO-Flow is \emph{consistent regardless of the value of $\alpha$}. 
As previously mentioned, this is expected as JKO-Flow is a proximal point algorithm that converges regardless of the step size $\alpha$.
In this case, five JKO-Flow iterations are enough to obtain this consistency.
Additional plots and hyperparameter setups for different benchmark datasets with similar performance results are shown in the Appendix. {Tab.~\ref{tab:synthetic_subproblem_alpha_comparison}} summarizes the comparison between the single shot and JKO-Flow on all synthetic 2D datasets for different values of $\alpha$.
We also show an illustration of all the datasets, estimated densities, and generated samples with JKO-Flow in Fig.~\ref{fig:toyModels}.

\begin{figure}[t] 
    \centering
    \begin{tabular}{ccccc} 	    
        $m = 3$ & $m = 4$ & $m = 5$  & $m = 8$  & $m = 16$ \\ 
        \\
        \multicolumn{5}{c}{\textbf{MMD}$^{\boldsymbol{2}}$ \textbf{Values, Single Shot}}
        \\
         1.1e-2 & 5.6e-3 & 2.46e-3 & 3.03e-3 & 2.7e-3 \\
        \includegraphics[width=0.17\textwidth]{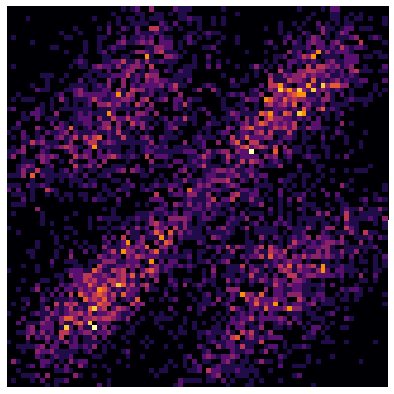}&
        \includegraphics[width=0.17\textwidth]{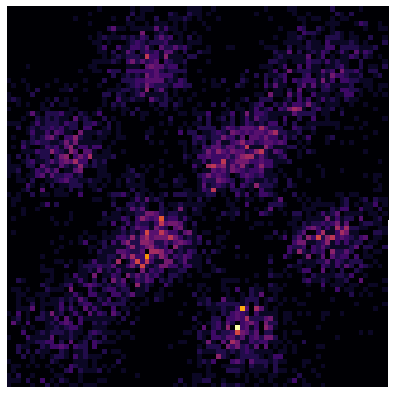}&
        \includegraphics[width=0.17\textwidth]{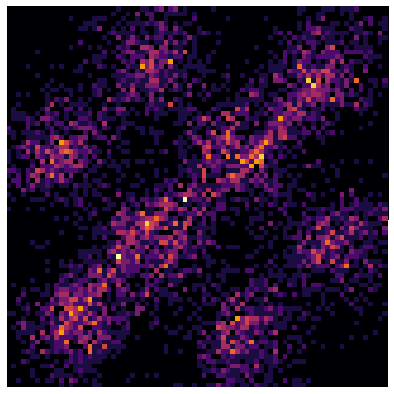}&
        \includegraphics[width=0.17\textwidth]{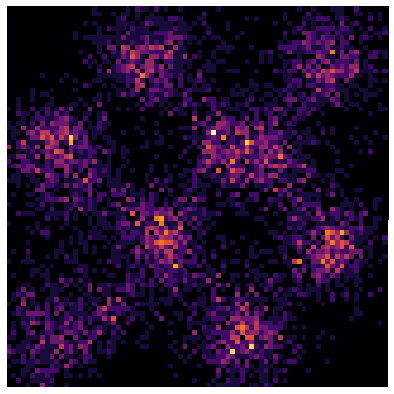}&
        \includegraphics[width=0.17\textwidth]{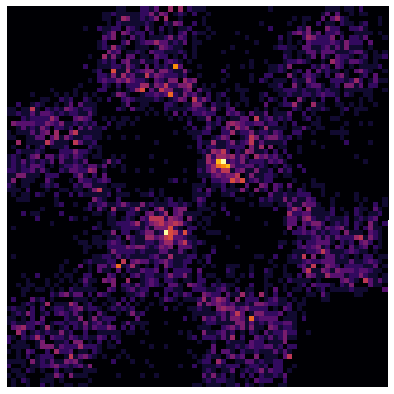}\\
        \\
        \multicolumn{5}{c}{\textbf{MMD}$^{\boldsymbol{2}}$ \textbf{Values, JKO-Flow (five iterations)}}
        \\
        5.6e-3 & 1.07e-3 & 2.7e-4 & 2.32e-4 & 4.16e-4 \\
        \includegraphics[width=0.17\textwidth]{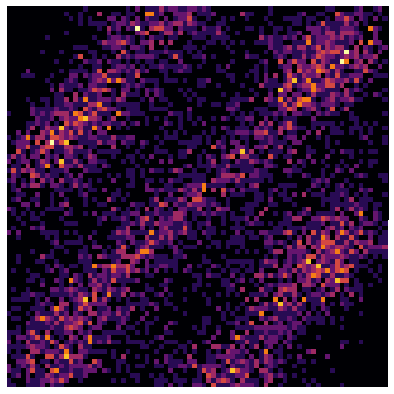}&
        \includegraphics[width=0.17\textwidth]{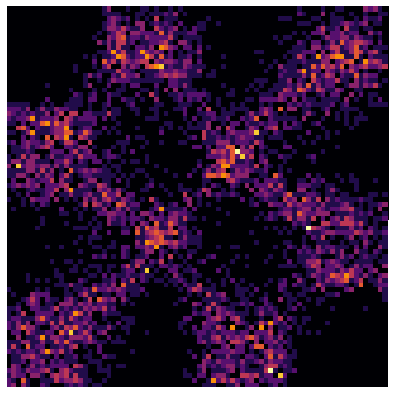}&
        \includegraphics[width=0.17\textwidth]{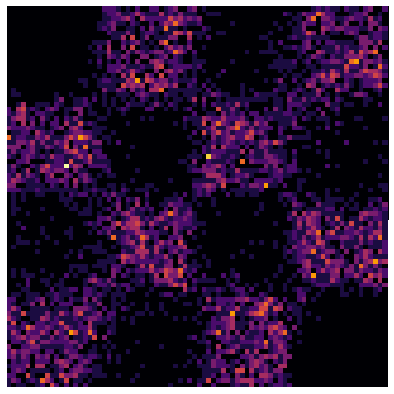}&
        \includegraphics[width=0.17\textwidth]{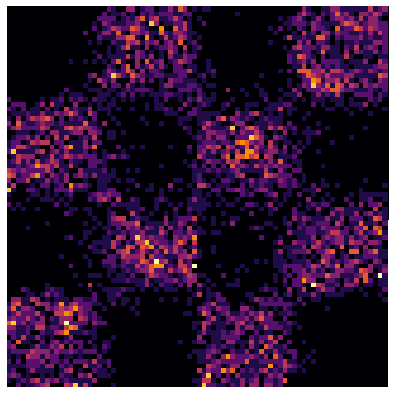}&
        \includegraphics[width=0.17\textwidth]{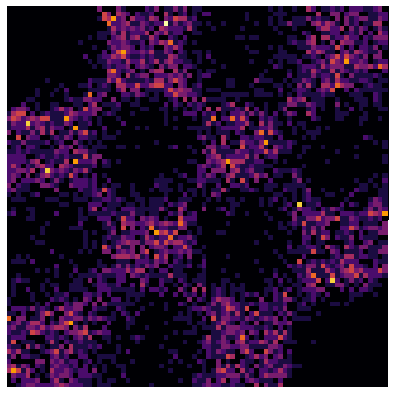} \\
    \end{tabular} 
    \centering 
    \\[2mm]
    \caption{Checkerboard dataset: Generated samples of $\widehat{\rho}_0$ using the standard single shot approach (top row). Generated samples using our proposed JKO-Flow using five iterations (bottom row). Here, we fix $\alpha = 5$ and vary the network width $m = 3, 4, 5, 8$, and $16$. JKO-Flow performs competitively even with fewer parameters.}
    \label{fig:checkerboardWidth}
\end{figure}

\noindent \textbf{Varying Network Size.}
In addition to obtaining consistent results for different values of $\alpha$, we also \emph{empirically} observe that JKO-Flow outperforms the single shot approach for different numbers of network parameters, i.e., network size. 
We illustrate this in Fig.~\ref{fig:checkerboardWidth}. 
This is also intuitive as we reformulate the problem of finding a single ``difficult'' optimal transportation problem as a sequence of ``smaller and easier'' OT problems. 
In this setup, we vary the width of a two-layer ResNet~\cite{he2016deep}. In particular, we choose the widths to be $m=3, 4, 5, 8$, and $16$. These correspond to $40, 53, 68, 125$, and $365$ parameters.  The hyperparameter $\alpha$ is chosen to be the best performing value for each synthetic dataset. 
All datasets vary $m$ for fixed $\alpha = 5$, except the 2 Spiral dataset, which uses $\alpha =50$; we chose these $\alpha$ values as they performed the best in the fixed $m$ experiments. 
Similar results are also shown for the remaining synthetic datasets in the appendix.
{Tab.~\ref{tab:synthetic_subproblem_networksize}} summarizes the comparison between the single shot and JKO-Flow on all synthetic 2D datasets.

\begin{figure}[t]
    \centering
    \addtolength{\tabcolsep}{-6pt} 
    \subfloat[\miniboone{} dimension 16 vs 17
    \label{fig:miniboone5.1}]{
	   	\begin{tabular}{ccc}
	   	
	   	{Samples} & Single Shot & \hspace{6pt} JKO-Flow \\
	   	$x \sim \rho_0(x)$ & $f(x)$ & \hspace{6pt}$f(x)$ \\
	   	\includegraphics[height=46pt]{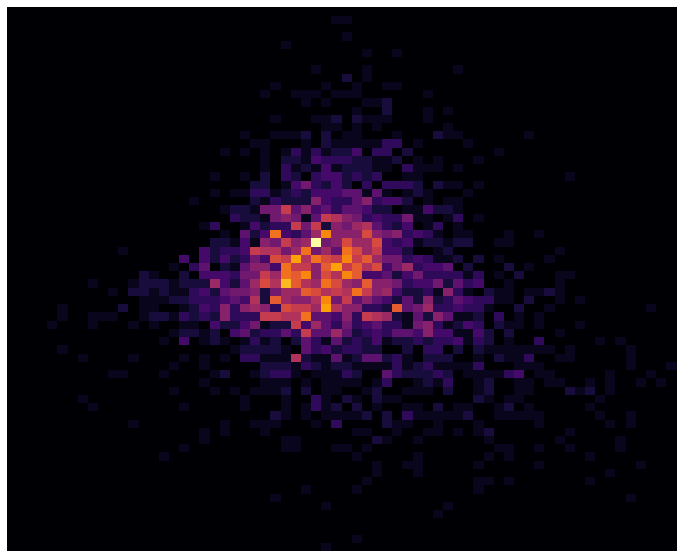}
	   	&
	   	\includegraphics[height=46pt]{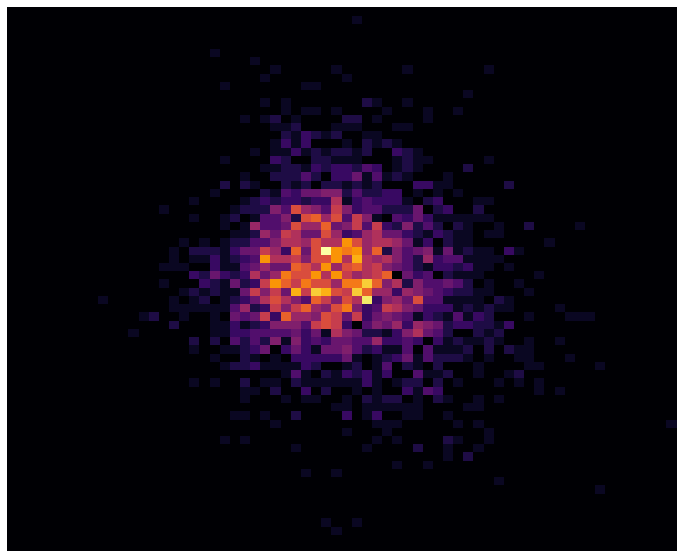}
	   	&
	   	\includegraphics[height=46pt]{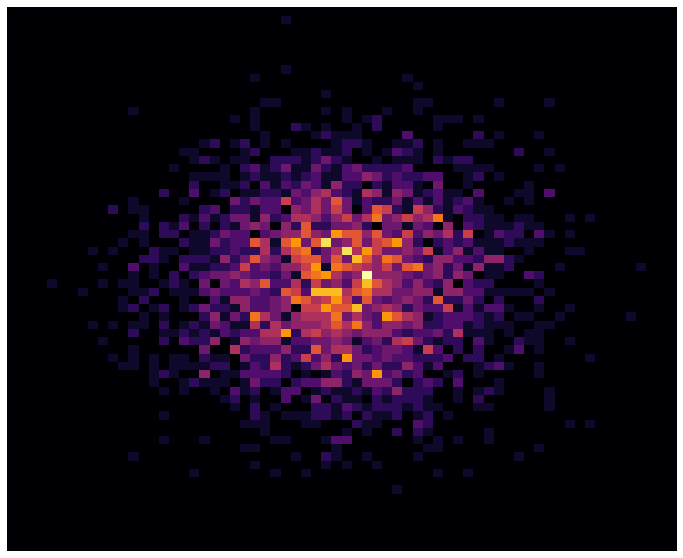} \\
	   	$y \sim \rho_1(y)$ & $f^{-1}(y)$ & \hspace{6pt} $f^{-1}(y)$ \\ 
	   	\includegraphics[height=46pt]{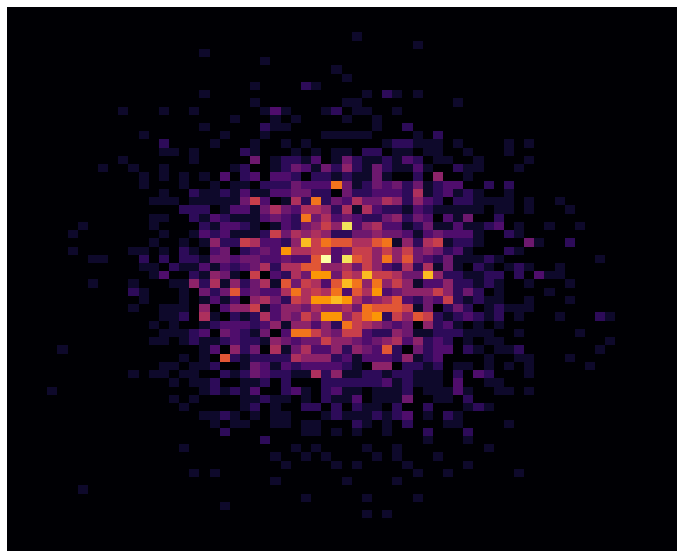}
	   	&
	   	\includegraphics[height=46pt]{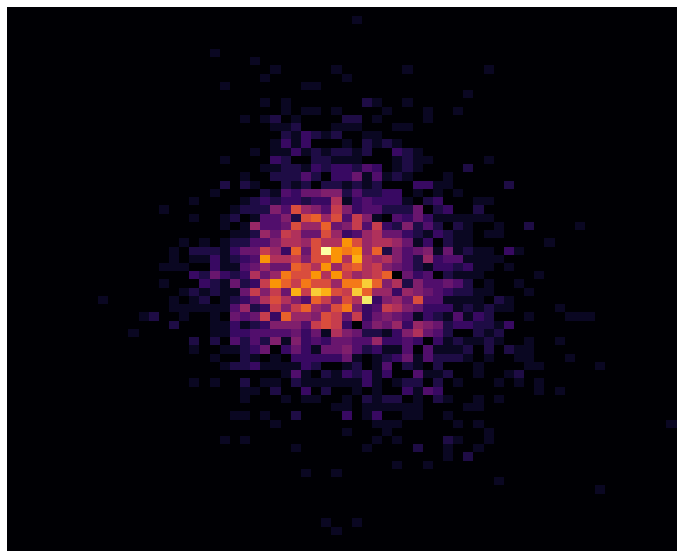}
	   	&
	   	\includegraphics[height=46pt]{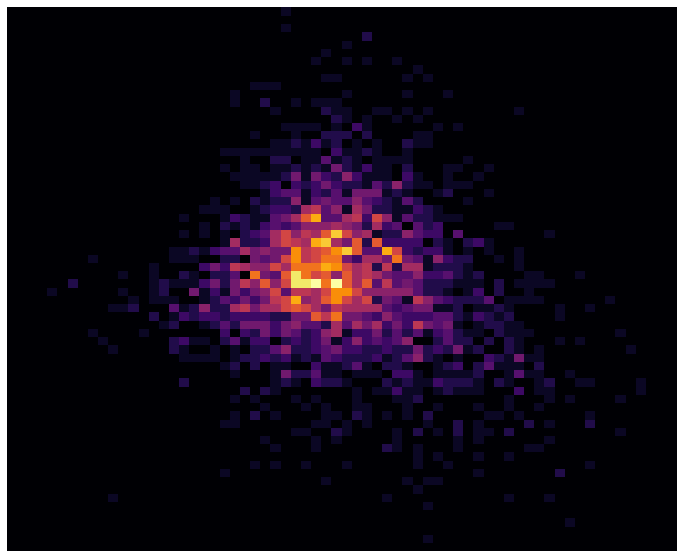}
	   	\end{tabular}%
	}\hspace{10pt}%
  	\subfloat[\miniboone{} dimension 28 vs 29
  	\label{fig:miniboone5.2}]{
	   	\begin{tabular}{ccc}
	   	{Samples} & Single Shot & \hspace{6pt} JKO-Flow \\
	   	$x \sim \rho_0(x)$ & $f(x)$ & \hspace{6pt} $f(x)$ \\
	   	\includegraphics[height=46pt]{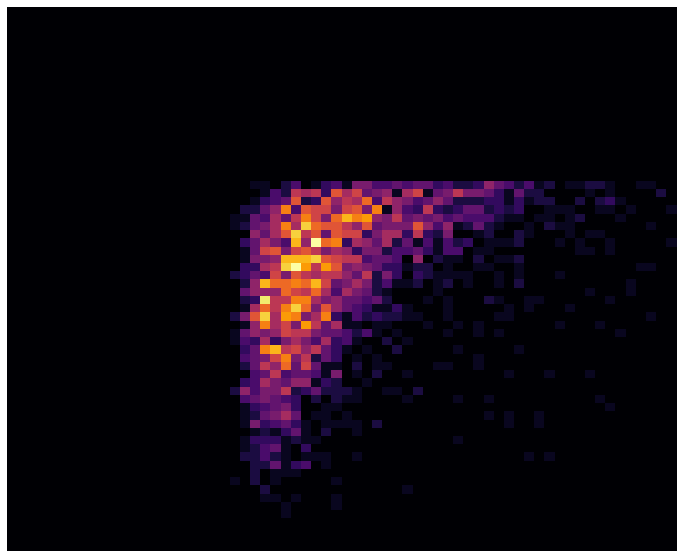}
	   	&
	   	\includegraphics[height=46pt]{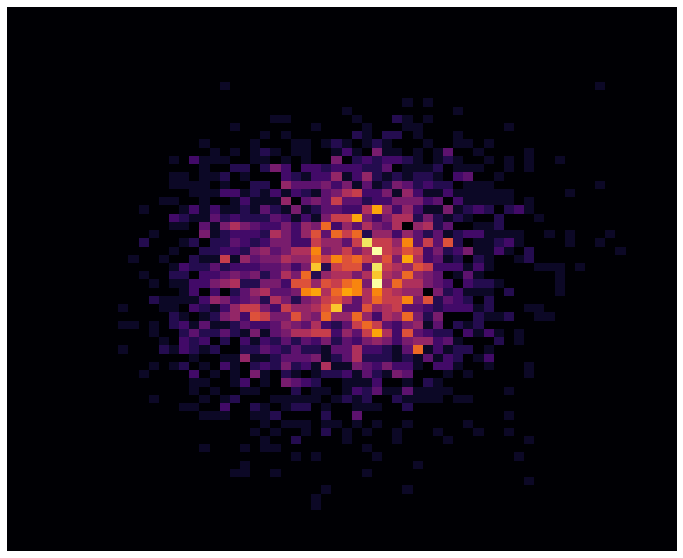}
	   	&
	   	\includegraphics[height=46pt]{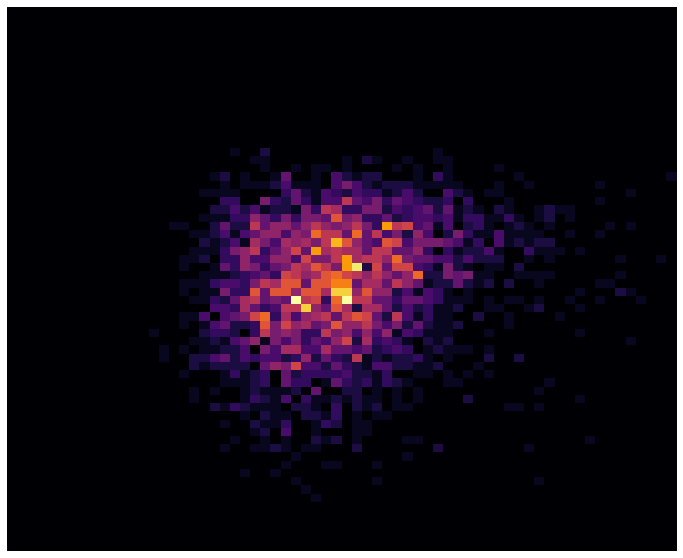} \\
	   	$y \sim \rho_1(y)$ & $f^{-1}(y)$ & \hspace{6pt} $f^{-1}(y)$ \\ 
	   	\includegraphics[height=46pt]{figures/miniboone_y.png}
	   	&
	   	\includegraphics[height=46pt]{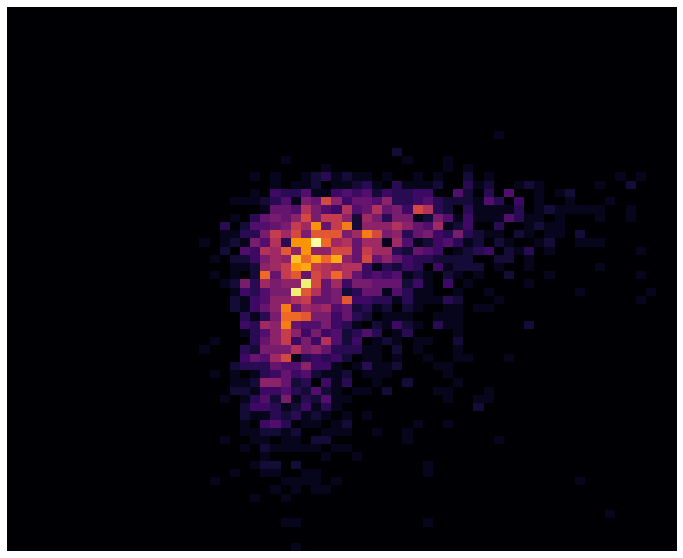}
	   	&
	   	\includegraphics[height=46pt]{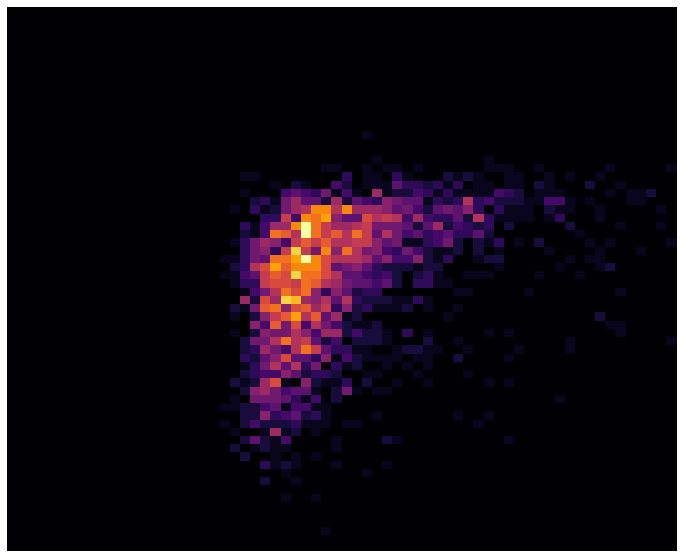}
	   	\end{tabular}%
	}%
    \caption{Generated samples for the 43-dimensional \miniboone{} dataset using the single shot approach and JKO-Flow with 10 iterations. To visualize the dataset, we show 2-dimensional slices. We show the forward flow $f(x)$ where $x \sim \rho_0$ and the genereated samples $f^{-1}(y)$ where $y \sim \rho_1$.
    }
    \label{fig:miniboone5}
\addtolength{\tabcolsep}{5pt} 
\end{figure}

\begin{figure}[t]
\centering
    \begin{tabular}{c c}
         \hspace{34pt} Dimension 16 vs. 17  & \hspace{34pt} Dimension 28 vs. 29 \\
            \includegraphics[width=0.47\linewidth]{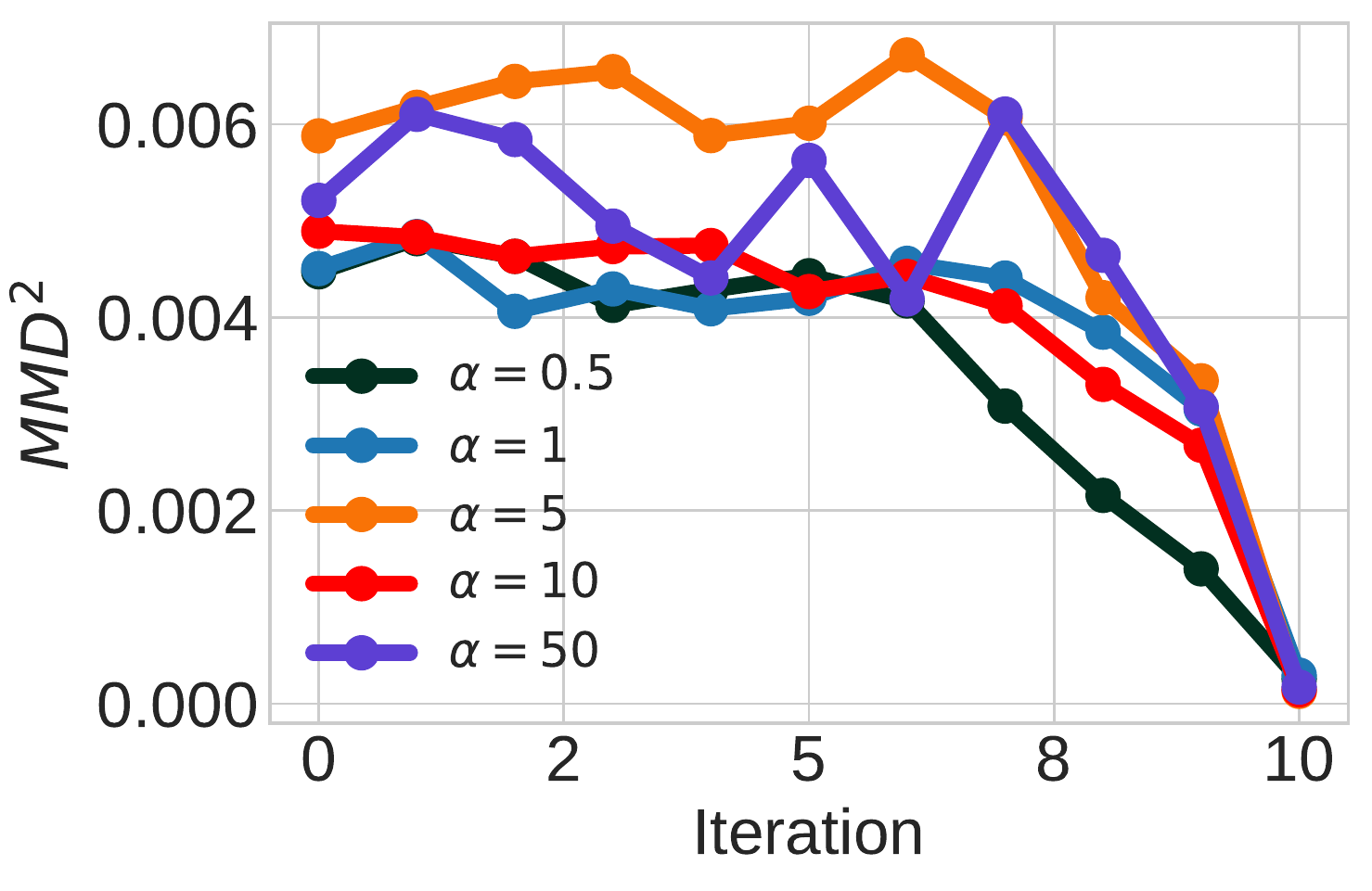} & \includegraphics[width=0.47\linewidth]{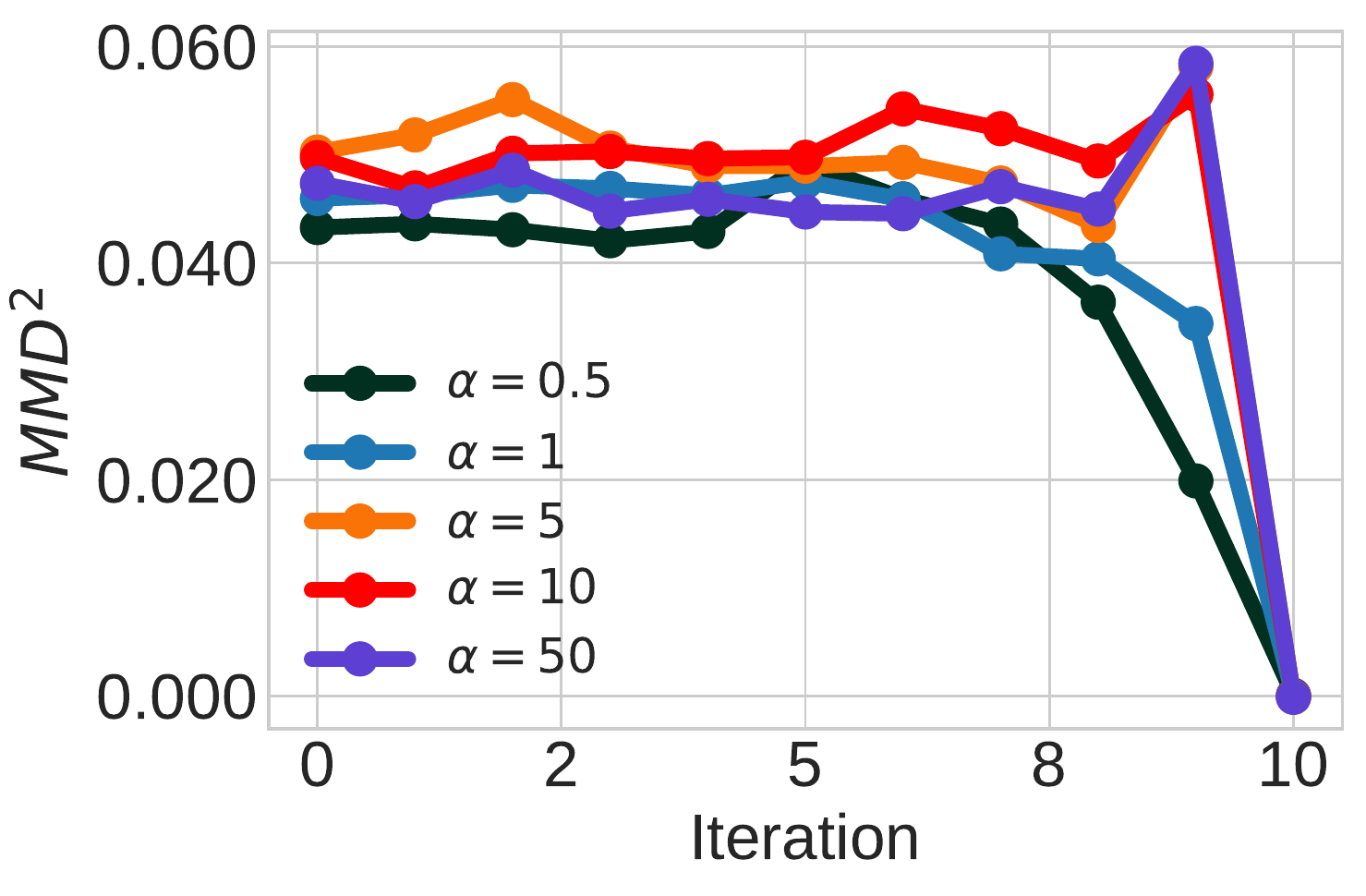}
    \end{tabular}
\caption{MMD$^2$ per iteration when using JKO-Flow after training. MMD$^2$ values vs. JKO-Flow iteration for 2-dimensional slice (dimensions 16-17 on the left and dimensions 28-29 on the right) of the \miniboone{} dataset. JKO-Flow achieves same accuracy \emph{regardless} of the value of $\alpha$.}
\label{fig: convergence}
\end{figure}

\noindent \textbf{Density Estimation on a Physics Dataset}
We train JKO-Flow on the 43-dimensional \miniboone{} dataset which is a high-dimensional, real-world physics dataset used as benchmark for high-dimensional density estimation algorithms in physics~\cite{misc_miniboone_particle_identification_199}. 
For this physics problem, our method is trained for $\alpha = 0.5, \, 1, \, 5, \, 10,\,50$ and using $10$ JKO-Flow iterations.
Fig~\ref{fig:miniboone5} shows generated samples with JKO-Flow and the standard single-shot approach for $\alpha = 5$. Since \miniboone{} is a high-dimensional dataset, we follow~\cite{onken2021ot} and plot two-dimensional slices. JKO-Flow generates better quality samples. Similar experiments for $\alpha=1, 10$, and $50$ are shown in Figs~\ref{fig:miniboone1},~\ref{fig:miniboone10}, and~\ref{fig:miniboone50} in the Appendix. Tab.~\ref{tab: miniboone_subproblem} summarizes the results for all values of $\alpha$. Note that we compute MMD values for all the dimensions as well as 2D slices; this is because we only have limited data (~3000 testing samples) and the 2D slice MMD give a better indication on the improvement of the generated samples.
Results show that \emph{the MMD is consistent across all $\alpha$ values for JKO-Flow}.
We also show the convergence (in MMD$^2$) of the miniboone dataset across each 2D slice in Fig.~\ref{fig: convergence}. As expected, smaller step size $\alpha$ values converge slower (see $\alpha = 0.5)$, but all converge to similar accuracy (unlike the single-shot).

\section{Conclusion}
We propose a new approach we call JKO-Flow to train OT-regularized CNFs without having to tune the regularization parameter $\alpha$. The key idea is to embed an underlying OT-based CNF solver into a Wasserstein gradient flow framework, also known as the JKO scheme; this approach makes the regularization parameter act as a ``time'' variable. Thus, instead of tuning $\alpha$,  we repeatedly solve proximal updates for a fixed (time variable) $\alpha$. 
In our setting, we choose OT-Flow~\cite{onken2021ot}, which leverages exact trace estimation for fast CNF training.
Our numerical experiments show that JKO-Flow leads to improved performance over the traditional approach. Moreover, \emph{JKO-Flow achieves similar results regardless of the choice of $\alpha$}. 
We also empirically observe improved performance when varying the size of the neural network. 
Future work will investigate JKO-Flow on similar problems such as deep learning-based methods for optimal control~\cite{fleming06controlled,onken2022neural,onken2021neural} and mean field games~\cite{lin2021alternating, ruthotto2020machine, agrawal2022random}.

\section*{Acknowledgments}
LN and SO were partially funded by AFOSR MURI FA9550-18-502, ONR N00014-18-1-2527, N00014-18-20-1-2093 and N00014-20-1-2787.

\begin{table}[t]
\centering
\fontsize{3}{4}\selectfont
\caption{Synthetic 2D Data: JKO-Flow performance for different values of $\alpha$. JKO-Flow returns consistent performance for different $\alpha$.}
\resizebox{\textwidth}{!}{%
\begin{tabular}{cccccc}
\hline
\scalebox{0.75}{\textbf{$\alpha$}} &  & \textbf{1} & \textbf{5} & \textbf{10} & \textbf{50} \\ \hline
\textbf{Dataset}   & \textbf{Approach}  & \multicolumn{4}{c}{\textbf{MMD\scalebox{0.45}{$\boldsymbol{^2}$}} } \\ \hline
\multirow{2}{*}{\textbf{Checkerboard}}   &  single Shot   & 3.58e-2    & 3.56e-3  & 1.42e-3   & 1.26e-3    \\ 
  &  JKO-Flow (5 iters)    & 4.9e-4  & 5.67e-4   & 6.40e-4  & 9.00e-4     \\ \hline
  \multirow{2}{*}{\textbf{2 Spirals}}   &  Single Shot   & 7.21e-2    & 2.30e-2  & 1.84e-2   & 7.73e-4    \\ 
  &  JKO-Flow (5 iters)    & 2.10e-2    & 4.62e-4   & 1.02e-4   & 5.37e-5     \\ \hline
\multirow{2}{*}{\textbf{Swiss Roll}}   &  Single Shot   & 4.74e-3   & 7.33e-4   & 2.86e-4   & 7.03e-4    \\ 
   &  JKO-Flow (5 iters)   & 5.16e-4   & 8.3e-5   & 3.27e-5    & 6.07e-4     \\ \hline
\multirow{2}{*}{\textbf{8 Gaussians}}   &  Single Shot & 9.18e-3  & 2.69e-4   & 3.94e-4    & 7.10e-4     \\ 
  &  JKO-Flow (5 iters)    & 1.07e-4   & 4.13e-5   & 2.67e-4  & 7.27e-6     \\ \hline
\multirow{2}{*}{\textbf{Circles}}   &  Single Shot  & 9.84e-3  & 2.24e-4    & 6.51e-4    & 1.04e-4     \\ 
  &  JKO-Flow (5 iters)    & 9.49e-4   & 9.97e-6   & 2.38e-5    & 9.28e-5     \\ \hline
  \multirow{2}{*}{\textbf{Pinwheel}}   &  Single Shot  & 1.18e-2 & 1.8e-3  & 1.37e-3  & 2.2e-5    \\ 
  &  JKO-Flow (5 iters)  & 4.7e-4   & 2.63e-4   & 3.84e-4    &  4.30e-4    \\ \hline
    \multirow{2}{*}{\textbf{Moons}}   &  Single Shot   & 1.45e-3  & 2.05e-3    & 2.49e-4    & 2.42e-4     \\ 
  &  JKO-Flow (5 iters)  & 1.92e-4  & 4.65e-5  & 4.3e-5   & 1.08e-4 \\\hline
\end{tabular}%
}
\label{tab:synthetic_subproblem_alpha_comparison}
\end{table}

\begin{table}[t]
\centering
\fontsize{3}{4}\selectfont
\caption{Synthetic 2D Data: Network width comparison for 1 and 5 iterations given a fixed, best performing $\alpha$. JKO-Flow performs better than the single shot approach for different network sizes.}
[2mm]
\resizebox{\textwidth}{!}{%
\begin{tabular}{ccccccc}
\hline
\scalebox{0.75}{\textbf{$m$}} &  & \textbf{3} & \textbf{4} & \textbf{5} & \textbf{8} & \textbf{16} \\ \hline
\textbf{Dataset}   & \textbf{Approach}  & \multicolumn{5}{c}{\textbf{MMD\scalebox{0.45}{$\boldsymbol{^2}$}} } \\ \hline
\multirow{2}{*}{\textbf{Checkerboard}}   &  Single Shot   & 1.10e-2   & 5.60e-3  & 2.46e-3   & 3.03e-3 & 2.70e-3    \\ 
  &  JKO-Flow (5 iters)    & 5.60e-3   & 1.07e-3    & 2.7e-4  & 2.32e-4  & 4.16e-4   \\ \hline
  \multirow{2}{*}{\textbf{2 Spirals}}   &  Single Shot   & 5.98e-3   & 4.54e-3  & 5.47e-3   & 1.19e-3 & 3.96e-3    \\ 
  &  JKO-Flow (5 iters)    & 1.42e-3    & 1.49e-5   & 6.11e-4   & 3.93e-5  & 2.19e-3   \\ \hline
\multirow{2}{*}{\textbf{Swiss Roll}}   &  Single Shot   & 8.89e-3   & 7.71e-3   & 1.41e-3   & 1.37e-3 & 1.52e-3   \\ 
   &  JKO-Flow (5 iters)   & 1.49e-3   & 2.90e-4  & 6.13e-4 & 2.29e-4  & 8.40e-5  \\ \hline
\multirow{2}{*}{\textbf{8 Gaussians}}   &  Single Shot & 2.20e-3 & 1.05e-3   & 1.04e-3   & 2.3e-4   & 5.05e-4 \\ 
  &  JKO-Flow (5 iters)    & 1.33e-4  & 9.85e-4   & 2.40e-5  & 3.96e-4 & 1.07e-4     \\ \hline
\multirow{2}{*}{\textbf{Circles}}   &  Single Shot  & 2.06e-3  & 1.72e-3    & 1.37e-3    & 1.69e-3 & 1.34e-3     \\ 
  &  JKO-Flow (5 iters)    & 1.94e-3 & 3.24e-4  & 7.71e-4  & 5.9e-5 & 1.01e-4    \\ \hline
    \multirow{2}{*}{\textbf{Pinwheel}} &  Single Shot  & 1.10e-2 & 4.03e-3  & 2.27e-3  & 3.80e-3 & 5.43e-4    \\ 
  &  JKO-Flow (5 iters)  & 1.20e-3 & 8.23e-4 & 1.60e-3  &  7.00e-5 & 2.69e-4   \\ \hline
  \multirow{2}{*}{\textbf{Moons}}   &  Single Shot   & 4.98e-3  & 4.54e-3  & 5.47e-3 & 1.2e-3 & 3.96e-3  \\ 
  &  JKO-Flow (5 iters)  & 1.42e-3  & 1.5e-5  & 6.11e-4  & 3.90e-5 & 2.19e-3 \\\hline

\end{tabular}%
}
\label{tab:synthetic_subproblem_networksize}
\end{table}

\begin{table}[t]
\centering
\fontsize{3}{4}\selectfont
\caption{\miniboone{}: Comparison of Single Shot and JKO-Flow for different values of $\alpha$.}
\label{tab: miniboone_subproblem}
[2mm]
\resizebox{\textwidth}{!}{%
\begin{tabular}{ccccccc}
\hline
\scalebox{0.75}{\textbf{$\alpha$}} &  &\textbf{0.5}  & \textbf{1} & \textbf{5} & \textbf{10} & \textbf{50} \\ \hline
\textbf{Dimensions}   & \textbf{Approach}  & \multicolumn{5}{c}{\textbf{MMD}\scalebox{0.45}{$\boldsymbol{^2}$}} \\ \hline
\multirow{2}{*}{\textbf{2d: 16 vs. 17}}   &  Single Shot & 4.62e-3    & 4.50e-3    & 5.88e-3    & 4.89e-3    & 5.2e-3     \\ 
  &  JKO-Flow (10 iters)  & 3.42e-4   & 2.94e-4   & 1.27e-4   & 1.43e-4    & 1.71e-4     \\ \hline
\multirow{2}{*}{\textbf{2d: 28 vs. 29}}   &  Single Shot  & 4.33e-2   & 4.60e-2   & 5.02e-2   & 4.97e-2    & 4.74e-2     \\ 
   &  JKO-Flow (10 iters) & 8.02e-5    & 3.31e-5  & 4.43e-5  & 6.33e-5 & 9.97e-5     \\ \hline
\multirow{2}{*}{\textbf{Full}}   &  Single Shot  & 4.7e-2  & 4.51e-3   & 4.75e-3    & 4.21e-3    & 4.27e-3     \\ 
  &  JKO-Flow (10 iters)  & 4.72e-4   & 4.72e-4   & 4.71e-4    & 4.72e-4    & 4.72e-4     \\ \hline

\end{tabular}%
}
\end{table}

\clearpage
\bibliographystyle{abbrv}
\bibliography{references}
\appendix

\section{Appendix}

\begin{figure}[ht] 
	\begin{tabular}{cccc} 
		$\alpha = 1$ & $\alpha = 5$  & $\alpha = 10$  & $\alpha = 50$ \\ 
              \\
        \multicolumn{4}{c}{\textbf{MMD$\boldsymbol{^2}$ Values, Single Shot}}
          \\
		7.21e-2 & 2.30e-2 & 1.84e-2 & 7.73e-4 \\
        \includegraphics[width=0.22\textwidth]{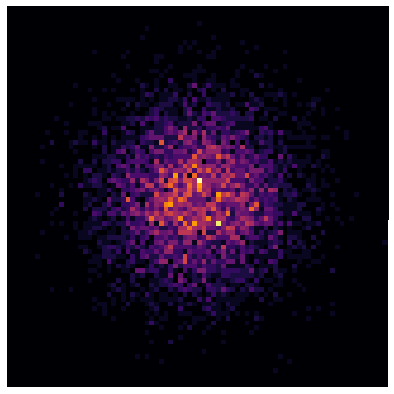}&
        \includegraphics[width=0.22\textwidth]{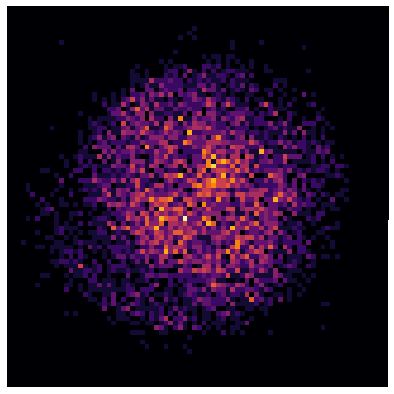}&
        \includegraphics[width=0.22\textwidth]{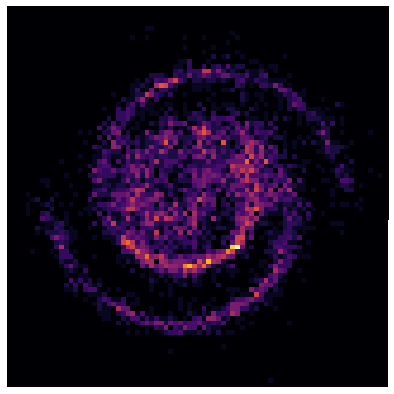}&
        \includegraphics[width=0.22\textwidth]{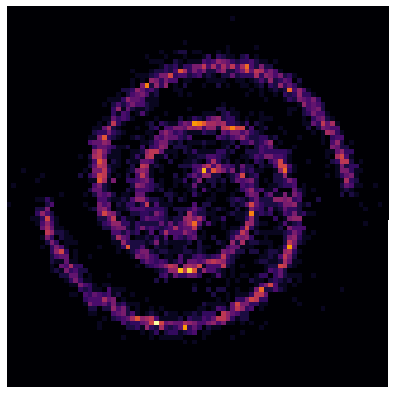}\\
                \\
        \multicolumn{4}{c}{\textbf{MMD$\boldsymbol{^2}$ Values, JKO-Flow (five iterations)}}
        \\
		2.10e-2 & 4.62e-4  & 1.02e-4 & 5.37e-5 \\
        \includegraphics[width=0.22\textwidth]{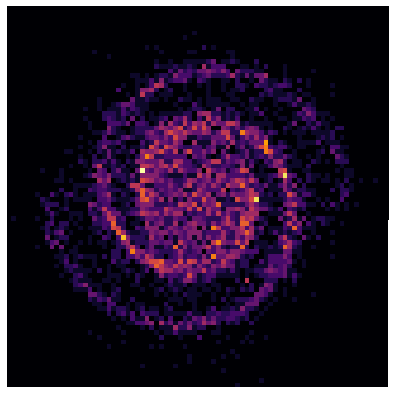}&
        \includegraphics[width=0.22\textwidth]{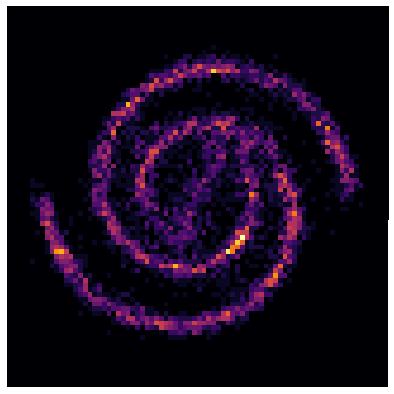}&
        \includegraphics[width=0.22\textwidth]{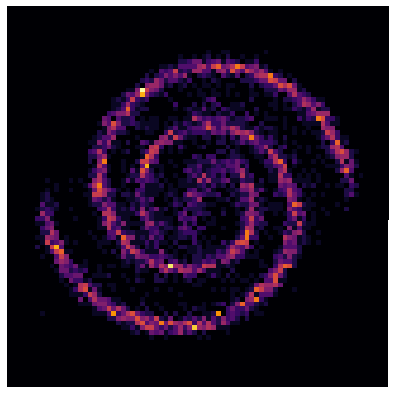}&
        \includegraphics[width=0.22\textwidth]{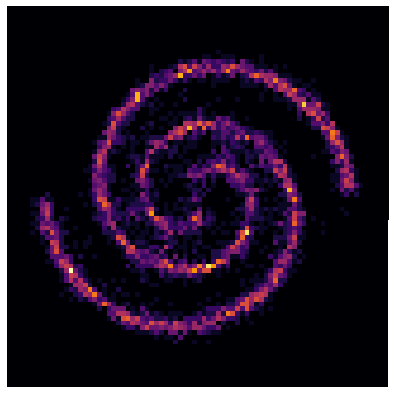} \\
	\end{tabular} 
	\centering 
	\\[2mm]
	\caption{2 Spirals dataset: Generated samples of $\hat{\rho}_0$ using the standard one-shot approach (top row). Generated using our proposed JKO-Flow using five iterations (bottom row). Here, we use $\alpha  = $ 1, 5, 10, 50. JKO-Flow returns consistent results \emph{regardless of the value of $\alpha$}.}
	
\end{figure} 

\begin{figure}[ht] 
	\begin{tabular}{cccc} 
		$\alpha = 1$ & $\alpha = 5$  & $\alpha = 10$  & $\alpha = 50$ \\ 
              \\
        \multicolumn{4}{c}{\textbf{MMD$\boldsymbol{^2}$ Values, Single Shot}}
          \\
		 4.73e-3 & 7.33e-4 & 2.86e-4 & 7.03e-4 \\
        \includegraphics[width=0.22\textwidth]{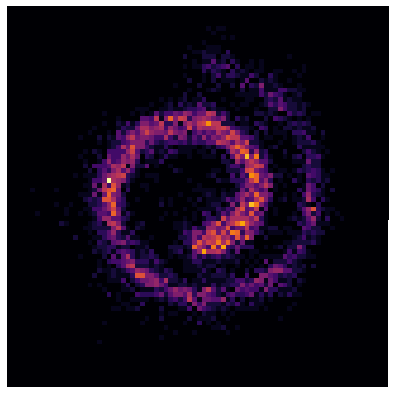}&
        \includegraphics[width=0.22\textwidth]{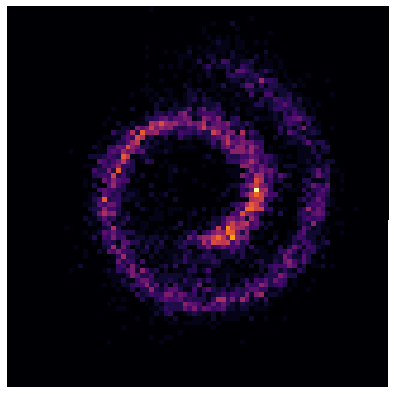}&
        \includegraphics[width=0.22\textwidth]{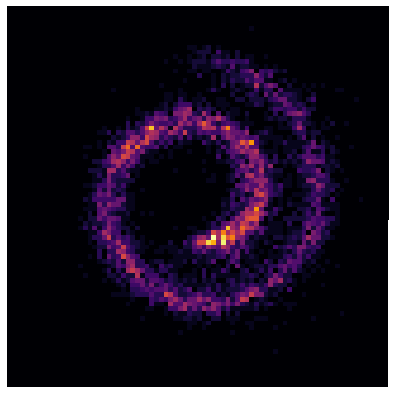}&
        \includegraphics[width=0.22\textwidth]{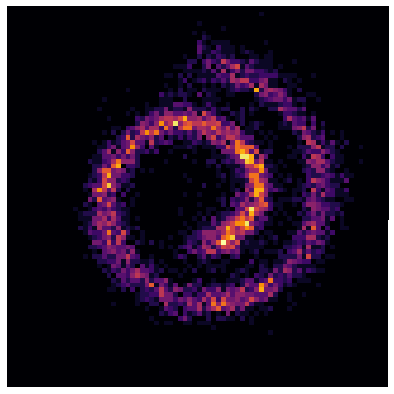}\\
                \\
        \multicolumn{4}{c}{\textbf{MMD$\boldsymbol{^2}$ Values, JKO-Flow (five iterations)}}
        \\
		5.16e-4 & 8.3e-5 & 3.27e-5 & 6.07e-4 \\
        \includegraphics[width=0.22\textwidth]{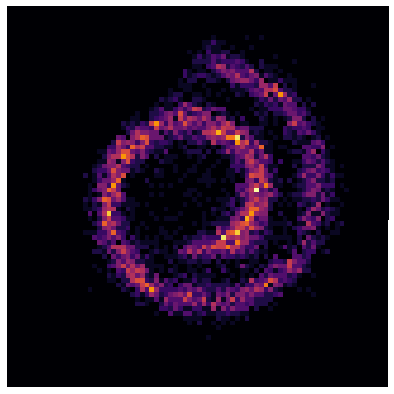}&
        \includegraphics[width=0.22\textwidth]{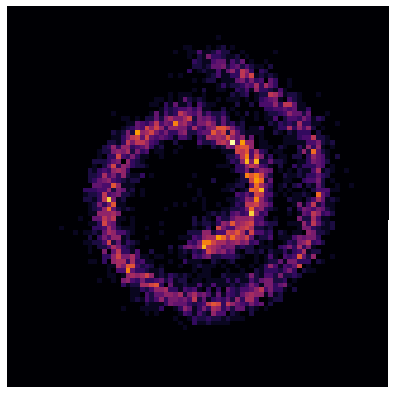}&
        \includegraphics[width=0.22\textwidth]{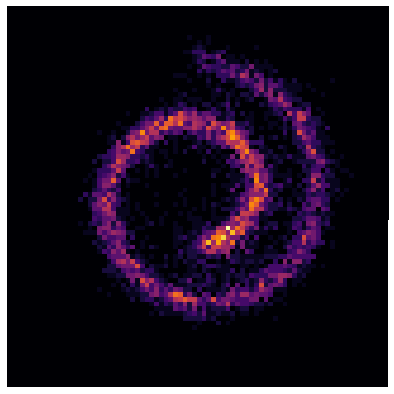}&
        \includegraphics[width=0.22\textwidth]{figures/spiral_1subproblems_alpha50_flow1.png} \\
	\end{tabular} 
	\centering 
	\\[2mm]
	\caption{Swiss Roll dataset: Generated samples of $\hat{\rho}_0$ using the standard one-shot approach (top row). Generated using our proposed JKO-Flow using five iterations (bottom row). Here, we use $\alpha  = $ 1, 5, 10, 50. JKO-Flow returns consistent results \emph{regardless of the value of $\alpha$}.}
	
\end{figure} 
\begin{figure}[ht] 
	\begin{tabular}{cccc} 
		$\alpha = 1$ & $\alpha = 5$  & $\alpha = 10$  & $\alpha = 50$ \\ 
  		\\ 
        \multicolumn{4}{c}{\textbf{MMD$\boldsymbol{^2}$ Values, Single Shot}}
		\\
		 9.18e-3 & 2.69e-4 & 3.94e-4 & 7.10e-4 \\
        \includegraphics[width=0.22\textwidth]{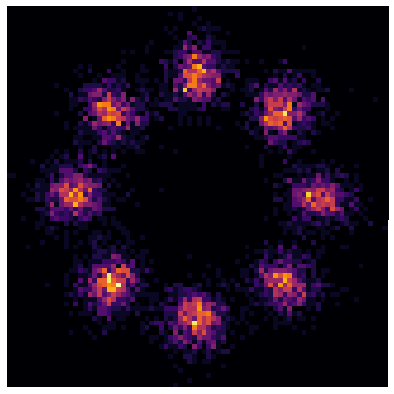}&
        \includegraphics[width=0.22\textwidth]{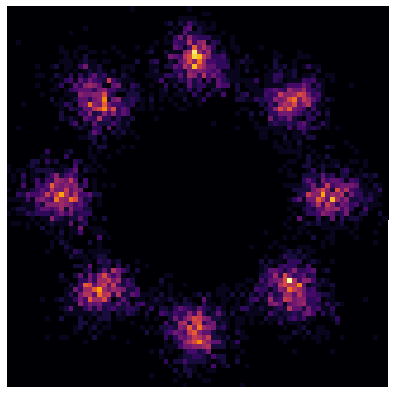}&
        \includegraphics[width=0.22\textwidth]{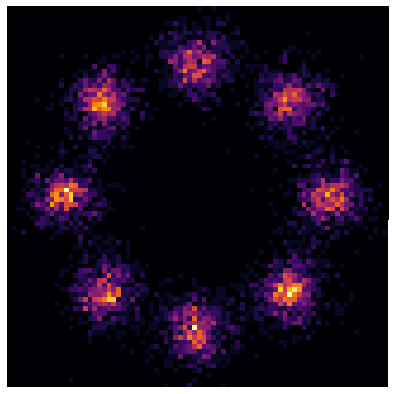}&
        \includegraphics[width=0.22\textwidth]{figures/8gaussians_1subproblems_alpha5_flow1.png}\
        \\
        \multicolumn{4}{c}{\textbf{MMD$\boldsymbol{^2}$ Values, JKO-Flow (five iterations)}}
        \\
		1.07e-4 & 4.13e-4 & 2.67e-4 & 7.27e-6 \\
        \includegraphics[width=0.22\textwidth]{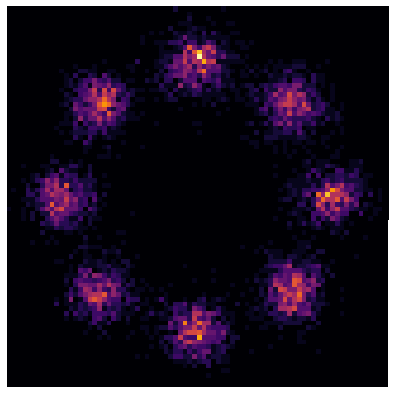}&
        \includegraphics[width=0.22\textwidth]{figures/8gaussians_1subproblems_alpha5_flow1.png}&
        \includegraphics[width=0.22\textwidth]{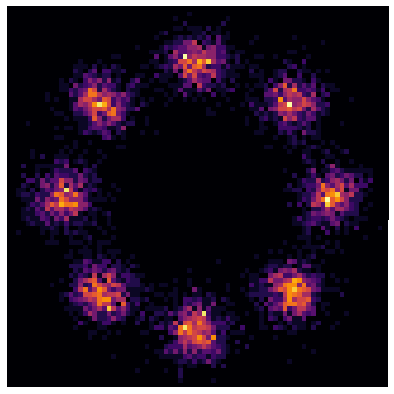}&
        \includegraphics[width=0.22\textwidth]{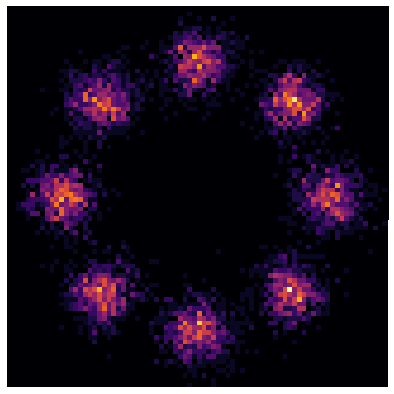} \\
	\end{tabular} 
	\centering 
	\\[2mm]
	\caption{8 Gaussians dataset: Generated samples of $\hat{\rho}_0$ using the standard one-shot approach (top row). Generated using our proposed JKO-Flow using five iterations (bottom row). Here, we use $\alpha  = $ 1, 5, 10, 50. JKO-Flow returns consistent results \emph{regardless of the value of $\alpha$}.}
	
\end{figure} 
\begin{figure}[ht] 
\centering
	\begin{tabular}{cccc} 
		$\alpha = 1$ & $\alpha = 5$  & $\alpha = 10$  & $\alpha = 50$ \\ 
      	\\ 
        \multicolumn{4}{c}{\textbf{MMD$\boldsymbol{^2}$ Values, Single Shot}}
		\\
		 9.84e-3 & 2.24e-4 & 6.51e-4 & 1.04e-4 \\
        \includegraphics[width=0.22\textwidth]{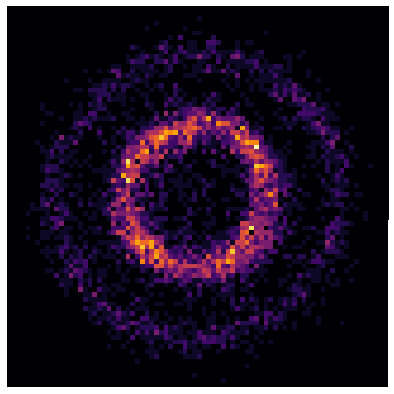}&
        \includegraphics[width=0.22\textwidth]{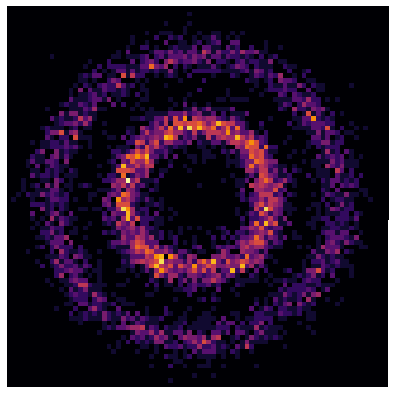}&
        \includegraphics[width=0.22\textwidth]{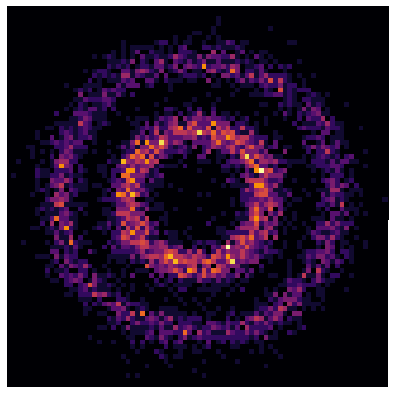}&
        \includegraphics[width=0.22\textwidth]{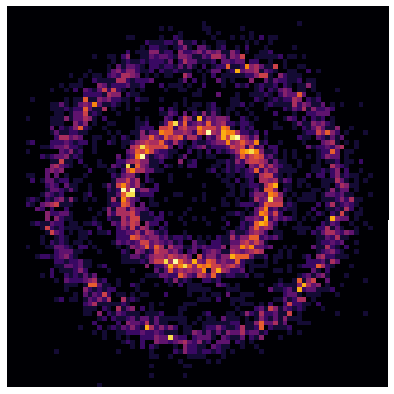}\\
        \\
        \multicolumn{4}{c}{\textbf{MMD$\boldsymbol{^2}$ Values, JKO-Flow (five iterations)}}
        \\
		9.49e-4 & 1.0e-5 & 2.4e-5 & 9.3e-5 \\
        \includegraphics[width=0.22\textwidth]{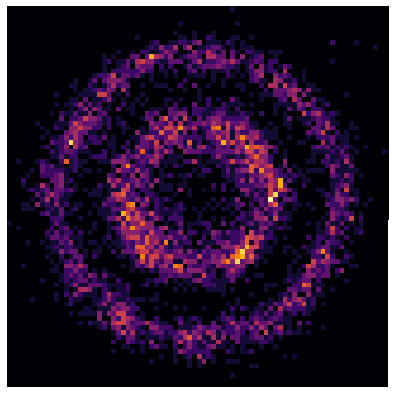}&
        \includegraphics[width=0.22\textwidth]{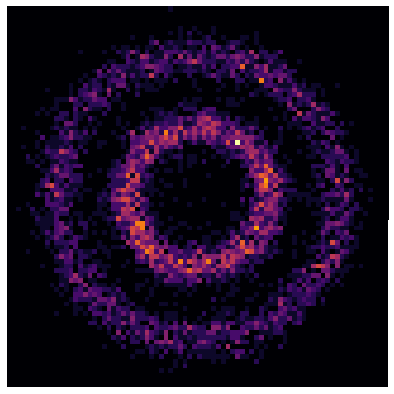}&
        \includegraphics[width=0.22\textwidth]{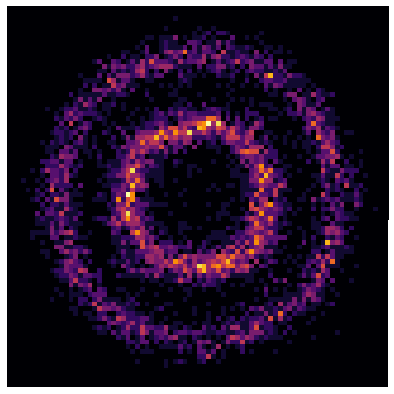}&
        \includegraphics[width=0.22\textwidth]{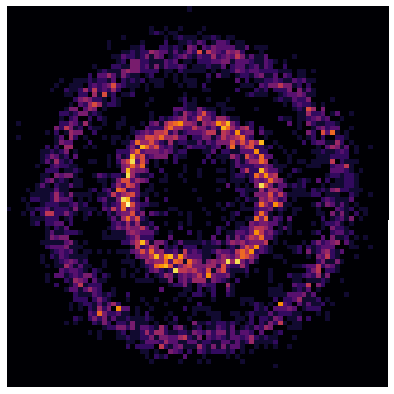} \\
	\end{tabular} 
	\centering 
	\\[2mm]
	\caption{Circles dataset: Generated samples of $\hat{\rho}_0$ using the standard one-shot approach (top row). Generated using our proposed JKO-Flow using five iterations (bottom row). Here, we use $\alpha  = $ 1, 5, 10, 50. JKO-Flow returns consistent results \emph{regardless of the value of $\alpha$}.}
	
\end{figure} 

\begin{table}[ht] 
\centering
	\begin{tabular}{cccc} 
		$\alpha = 1$ & $\alpha = 5$  & $\alpha = 10$  & $\alpha = 50$ \\ 
        \\ 
        \multicolumn{4}{c}{\textbf{MMD$\boldsymbol{^2}$ Values, Single Shot}}
		\\
		  1.18e-2 & 1.8e-3 & 1.37e-4 & 2.2e-5 \\
        \includegraphics[width=0.22\textwidth]{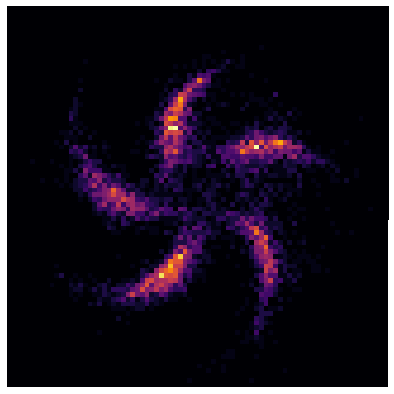}&
        \includegraphics[width=0.22\textwidth]{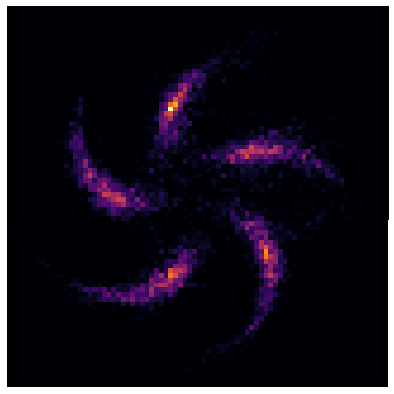}&
        \includegraphics[width=0.22\textwidth]{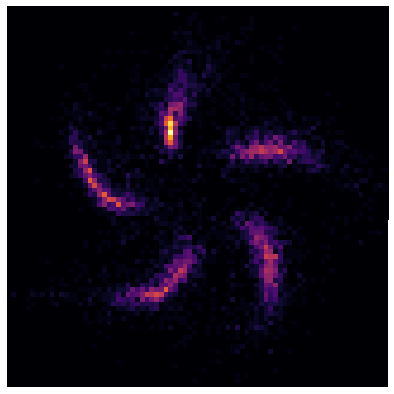}&
        \includegraphics[width=0.22\textwidth]{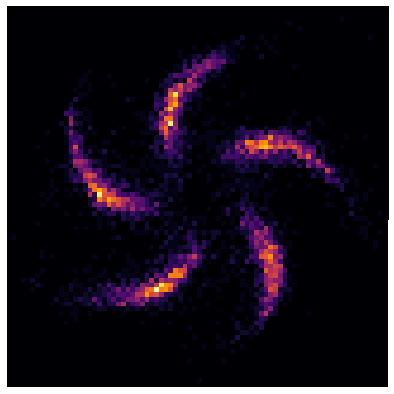}\\
        \\
        \multicolumn{4}{c}{\textbf{MMD$\boldsymbol{^2}$ Values, JKO-Flow (five iterations)}}
        \\
		4.78e-4 & 2.63e-4 & 3.84e-4 & 4.30e-4 \\
        \includegraphics[width=0.22\textwidth]{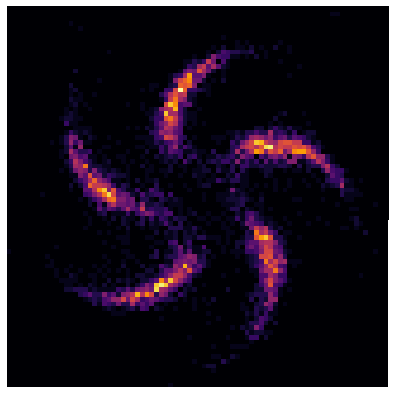}&
        \includegraphics[width=0.22\textwidth]{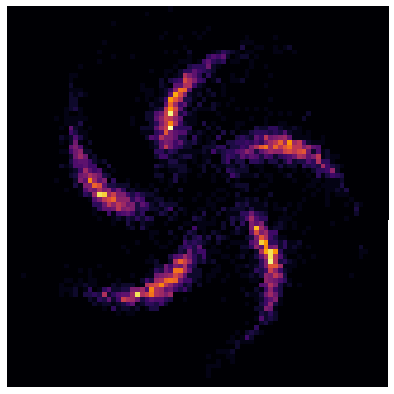}&
        \includegraphics[width=0.22\textwidth]{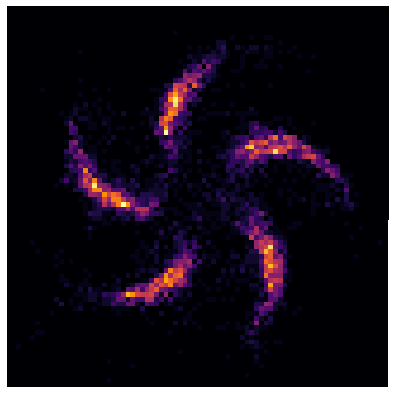}&
        \includegraphics[width=0.22\textwidth]{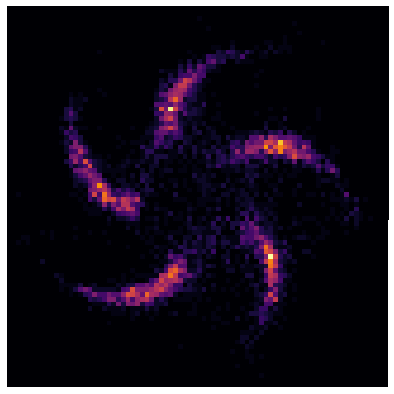} \\
	\end{tabular} 
	\centering 
	\\[2mm]
	\caption{Pinwheel dataset: Generated samples of $\hat{\rho}_0$ using the standard one-shot approach (top row). Generated using our proposed JKO-Flow using five iterations (bottom row). Here, we use $\alpha  = $ 1, 5, 10, 50. JKO-Flow returns consistent results \emph{regardless of the value of $\alpha$}.}
	
\end{table} 
\begin{figure}[ht] 
	\begin{tabular}{cccc} 
		$\alpha = 1$ & $\alpha = 5$  & $\alpha = 10$  & $\alpha = 50$ \\ 
        \\ 
        \multicolumn{4}{c}{\textbf{MMD$\boldsymbol{^2}$ Values, Single Shot}}
		\\
		  1.45e-3 & 2.05e-3 & 2.49e-4 & 2.42e-4 \\
        \includegraphics[width=0.22\textwidth]{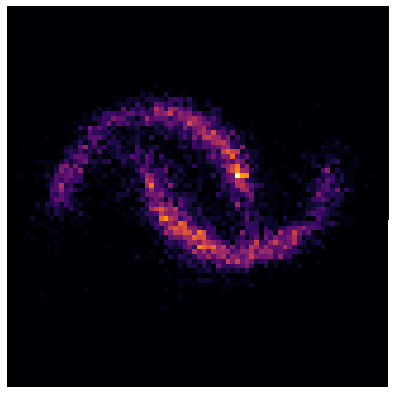}&
        \includegraphics[width=0.22\textwidth]{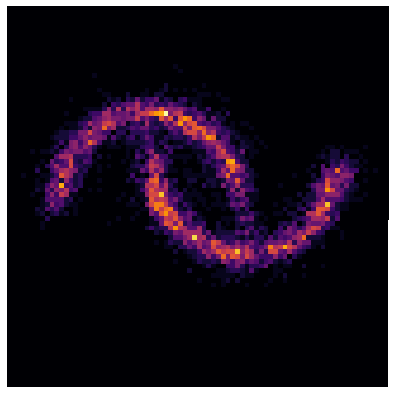}&
        \includegraphics[width=0.22\textwidth]{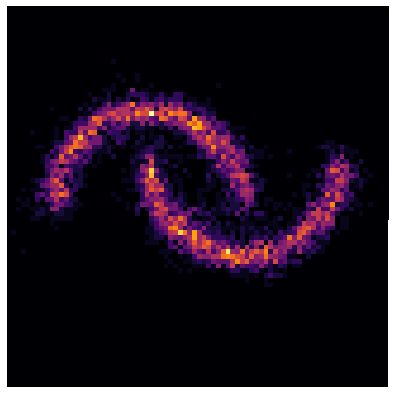}&
        \includegraphics[width=0.22\textwidth]{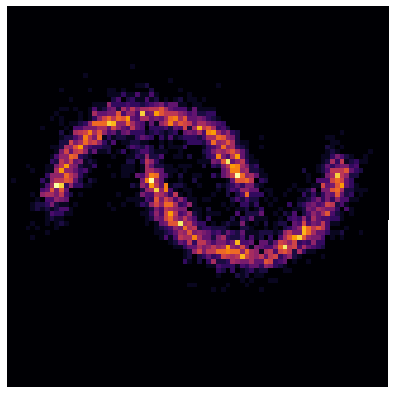}\\
        \\
        \multicolumn{4}{c}{\textbf{MMD$\boldsymbol{^2}$ Values, JKO-Flow (five iterations)}}
        \\
		1.9e-5 & 5.0e-5 & 4.3e-5 & 1.08e-4 \\
        \includegraphics[width=0.22\textwidth]{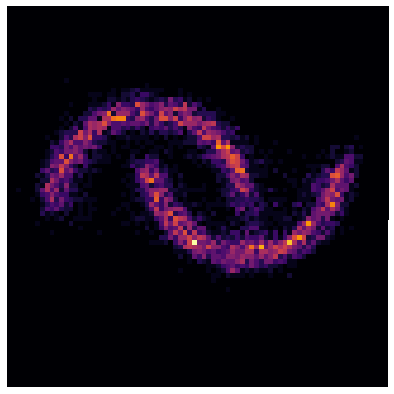}&
        \includegraphics[width=0.22\textwidth]{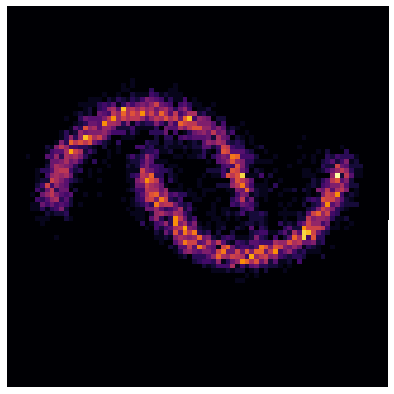}&
        \includegraphics[width=0.22\textwidth]{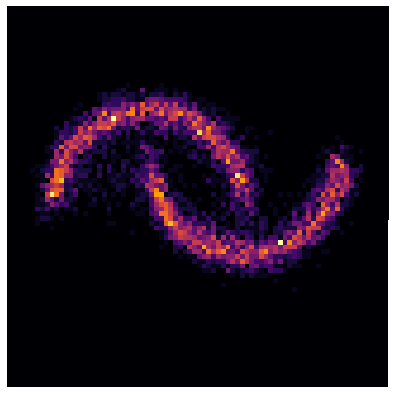}&
        \includegraphics[width=0.22\textwidth]{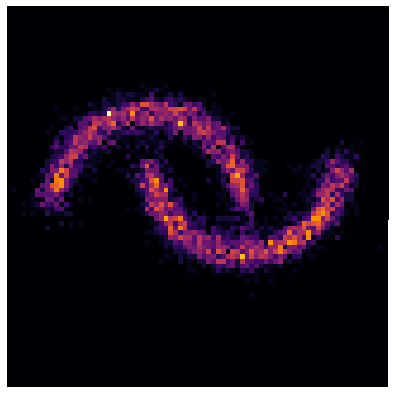} \\
	\end{tabular} 
	\centering 
	\\[2mm]
	\caption{Moons dataset: Generated samples of $\hat{\rho}_0$ using the standard one-shot approach (top row). Generated using our proposed JKO-Flow using five iterations (bottom row). Here, we use $\alpha  = $ 1, 5, 10, 50. JKO-Flow returns consistent results \emph{regardless of the value of $\alpha$}.}
\end{figure} 
\begin{figure}[ht] 
    \centering
	\begin{tabular}{ccccc} 	
		$m = 3$ & $m = 4$ & $m = 5$  & $m = 8$  & $m = 16$ 
		\\ 
        \multicolumn{5}{c}{\textbf{MMD$\boldsymbol{^2}$ Values, Single Shot}}
		\\
		  5.98e-3 & 4.54e-3 & 5.47e-3 & 1.19e-3 & 3.96e-3 \\
        \includegraphics[width=0.17\textwidth]{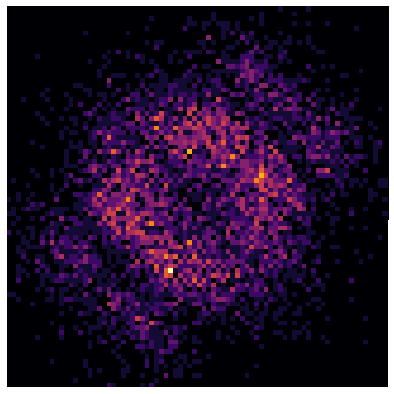}&
        \includegraphics[width=0.17\textwidth]{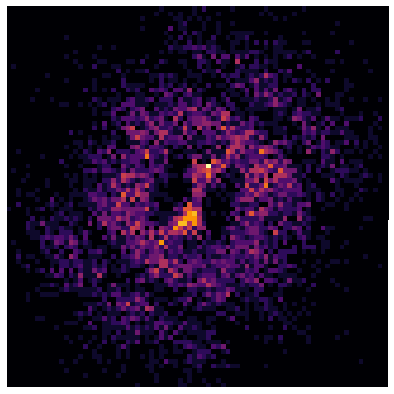}&
        \includegraphics[width=0.17\textwidth]{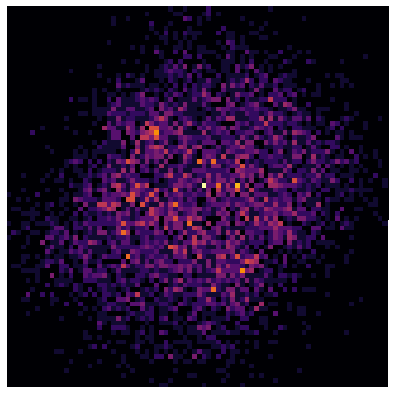}&
        \includegraphics[width=0.17\textwidth]{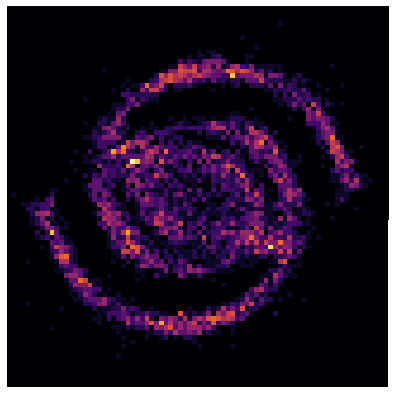}&
        \includegraphics[width=0.17\textwidth]{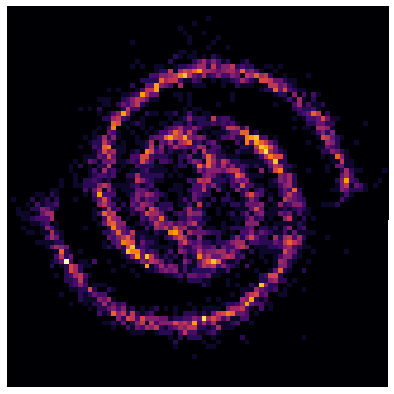}\\
         \\
        \multicolumn{5}{c}{\textbf{MMD$\boldsymbol{^2}$ Values, JKO-Flow (five iterations)}}
        \\
		1.42e-3 & 1.49e-5 & 6.11e-4 & 3.93e-5 & 2.19e-3 \\
        \includegraphics[width=0.17\textwidth]{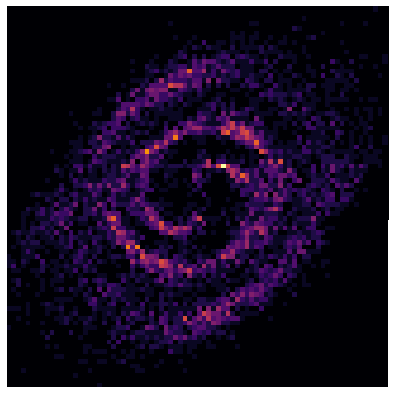}&
        \includegraphics[width=0.17\textwidth]{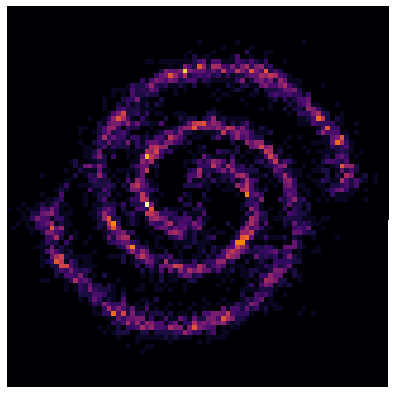}&
        \includegraphics[width=0.17\textwidth]{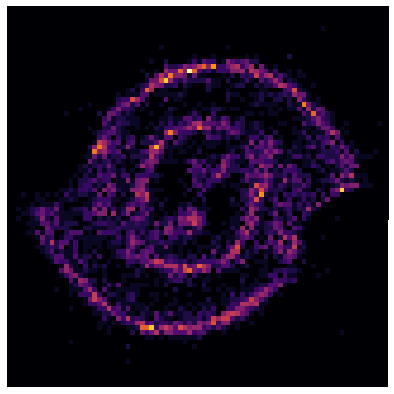}&
        \includegraphics[width=0.17\textwidth]{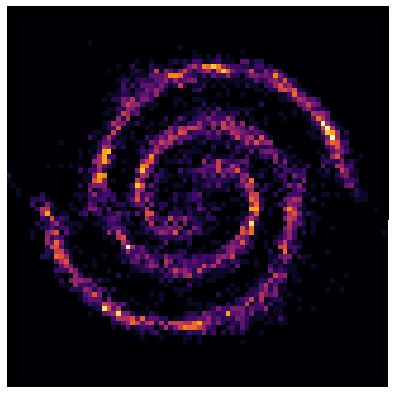}&
        \includegraphics[width=0.17\textwidth]{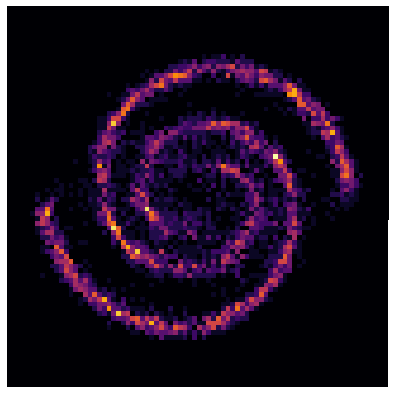} \\
	\end{tabular} 
	\centering 
	\\[2mm]
	\caption{2 Spirals dataset: Generated samples of $\hat{\rho}_0$ using the standard single shot approach (top row). Generated samples using our proposed JKO-Flow using five iterations (bottom row). Here, we fix $\alpha = 50$ and vary the network width $m = 3, 4, 5, 8$, and $16$. JKO-Flow performs competitively even with lower number of parameters.}
	
\end{figure} 
\begin{figure}[ht] 
    \centering
	\begin{tabular}{ccccc} 	
		$m = 3$ & $m = 4$ & $m = 5$  & $m = 8$  & $m = 16$ 
		\\ 
        \multicolumn{5}{c}{\textbf{MMD$\boldsymbol{^2}$ Values, Single Shot}}
		\\
		 8.89e-3 & 7.71e-3 & 1.41e-3 & 1.37e-3 & 1.52e-3 \\
        \includegraphics[width=0.17\textwidth]{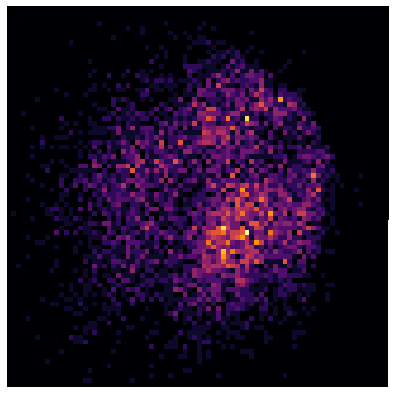}&
        \includegraphics[width=0.17\textwidth]{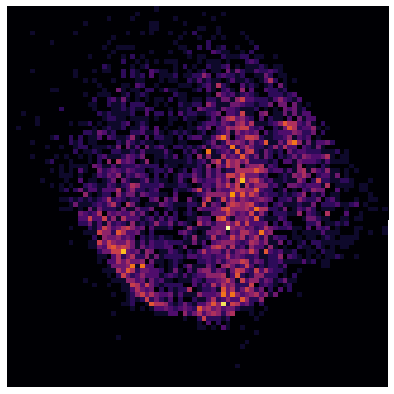}&
        \includegraphics[width=0.17\textwidth]{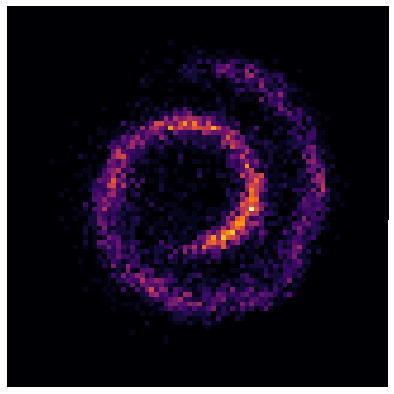}&
        \includegraphics[width=0.17\textwidth]{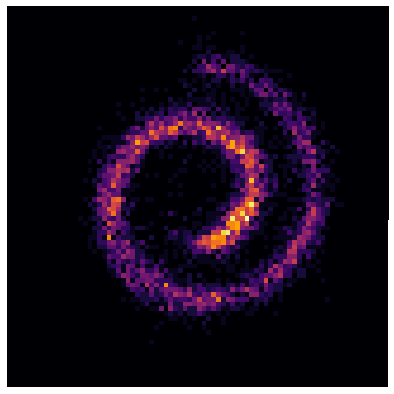}&
        \includegraphics[width=0.17\textwidth]{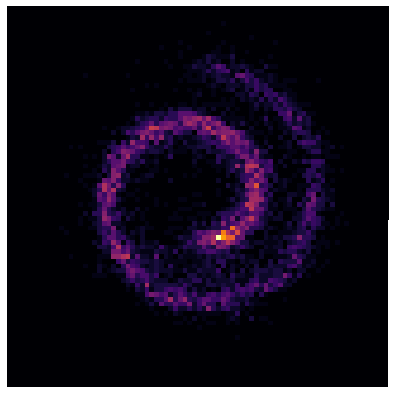}\\
         \\
        \multicolumn{5}{c}{\textbf{MMD$\boldsymbol{^2}$ Values, JKO-Flow (five iterations)}}
        \\
		1.49e-3 & 2.90e-4 & 6.13e-4 & 2.29e-4 & 8.4e-5 \\
        \includegraphics[width=0.17\textwidth]{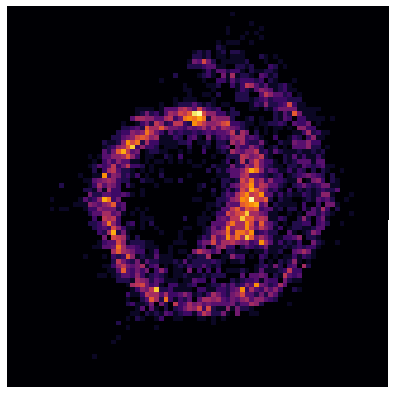}&
        \includegraphics[width=0.17\textwidth]{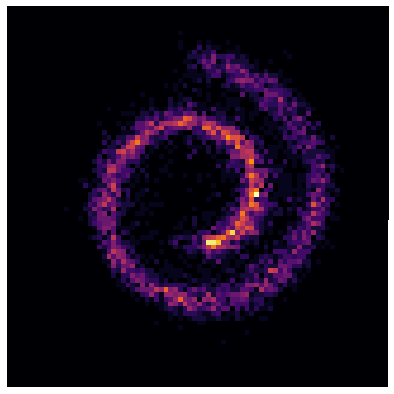}&
        \includegraphics[width=0.17\textwidth]{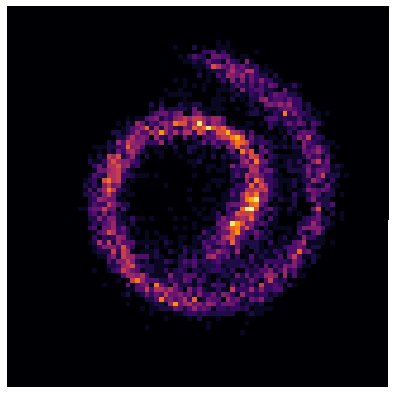}&
        \includegraphics[width=0.17\textwidth]{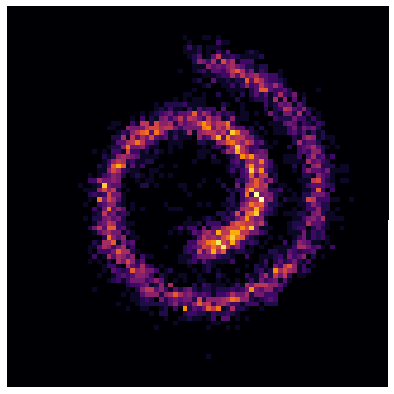}&
        \includegraphics[width=0.17\textwidth]{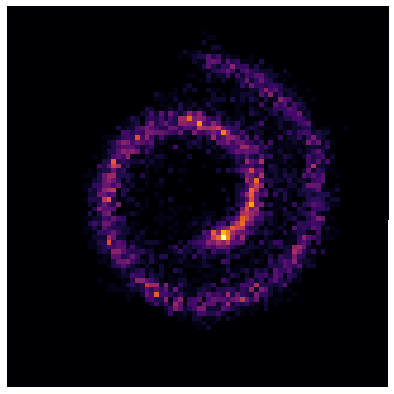} \\
	\end{tabular} 
	\centering 
	\\[2mm]
	\caption{Swiss Roll dataset: Generated samples of $\hat{\rho}_0$ using the standard single shot approach (top row). Generated samples using our proposed JKO-Flow using five iterations (bottom row). Here, we fix $\alpha = 5$ and vary the network width $m = 3, 4, 5, 8$, and $16$. JKO-Flow performs competitively even with lower number of parameters.}
	
\end{figure} 
\begin{figure}[ht] 
	\begin{tabular}{ccccc} 	    
		$m = 3$ & $m = 4$ & $m = 5$  & $m = 8$  & $m = 16$ \\ 
    	\\ 
        \multicolumn{5}{c}{\textbf{MMD$\boldsymbol{^2}$ Values, Single Shot}}
		\\
		  2.2e-3 & 1.05e-3 & 1.04e-3 &2.3e-4 & 5.05e-4 \\
		\includegraphics[width=0.17\textwidth]{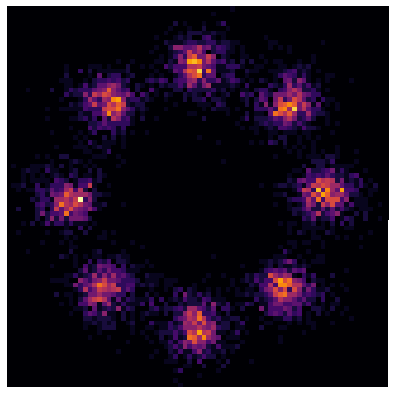}&
        \includegraphics[width=0.17\textwidth]{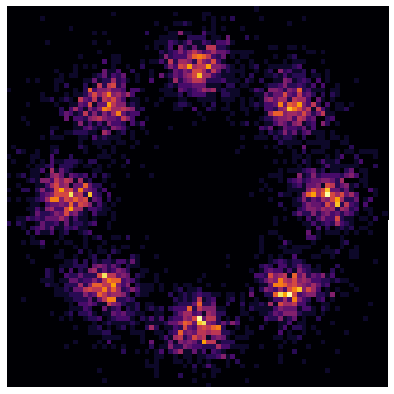}&
        \includegraphics[width=0.17\textwidth]{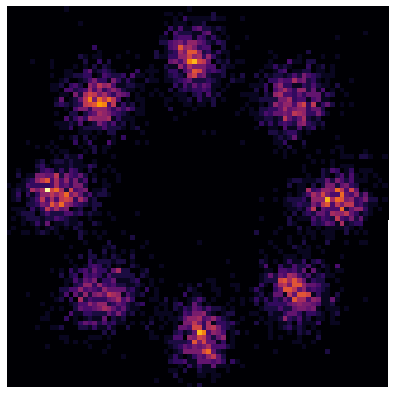}&
        \includegraphics[width=0.17\textwidth]{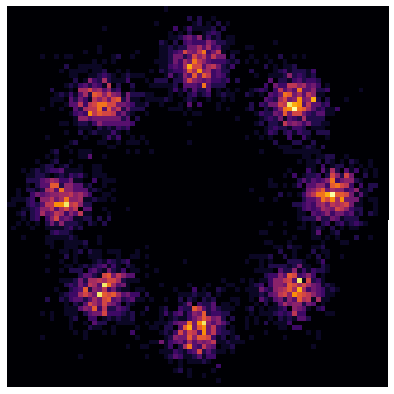}&
        \includegraphics[width=0.17\textwidth]{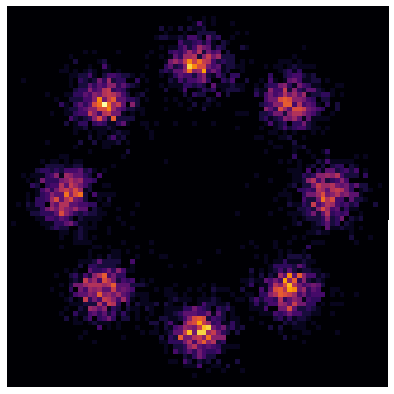}\\
        \\
        \multicolumn{5}{c}{\textbf{MMD$\boldsymbol{^2}$ Values, JKO-Flow (five iterations)}}
        \\
		1.33e-4 & 9.85e-4 & 2.4e-5 & 3.96e-4 & 1.07e-4 \\
		\includegraphics[width=0.17\textwidth]{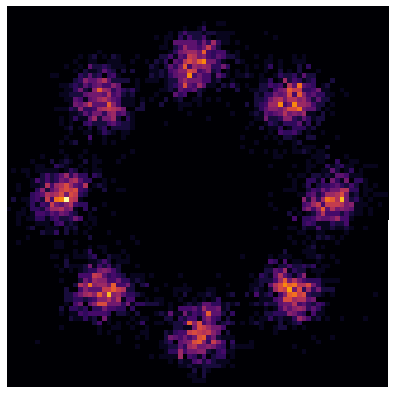}&
        \includegraphics[width=0.17\textwidth]{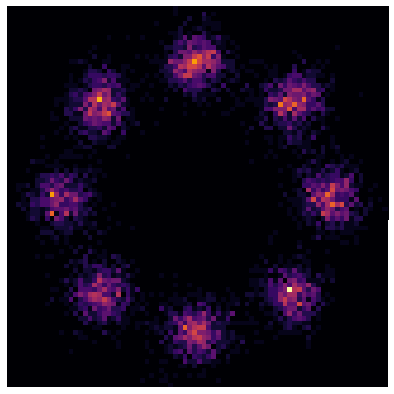}&
        \includegraphics[width=0.17\textwidth]{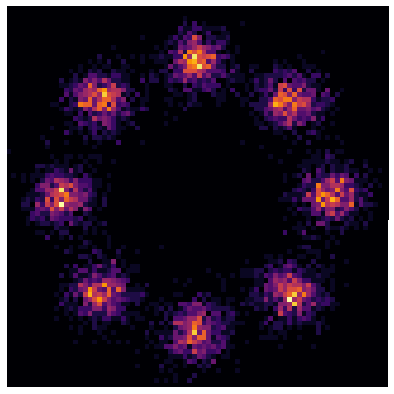}&
        \includegraphics[width=0.17\textwidth]{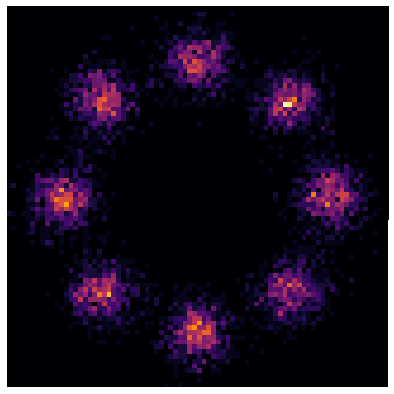}&
        \includegraphics[width=0.17\textwidth]{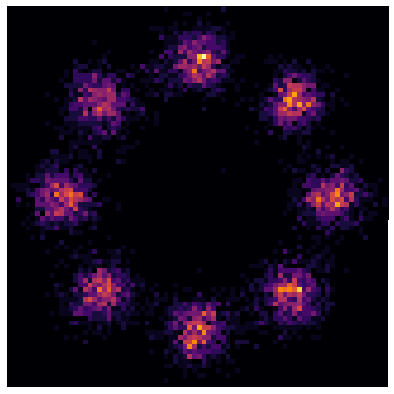} \\
	\end{tabular} 
	\centering 
	\\[2mm]
	\caption{8 Gaussians dataset: Generated samples of $\hat{\rho}_0$ using the standard single shot approach (top row). Generated samples using our proposed JKO-Flow using five iterations (bottom row). Here, we fix $\alpha = 5$ and vary the network width $m = 3, 4, 5, 8$, and $16$. JKO-Flow performs competitively even with lower number of parameters.}
	
\end{figure} 

\begin{figure}[ht] 
\centering
	\begin{tabular}{ccccc} 	    
		$m = 3$ & $m = 4$ & $m = 5$  & $m = 8$  & $m = 16$ \\ 
        \\ 
        \multicolumn{5}{c}{\textbf{MMD$\boldsymbol{^2}$ Values, Single Shot}}
		\\
		 2.06e-3 & 1.72e-3 & 1.37e-3 & 1.69e-3 & 1.34e-3 \\
		\includegraphics[width=0.17\textwidth]{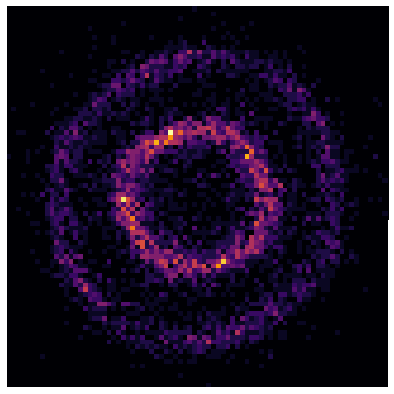}&
        \includegraphics[width=0.17\textwidth]{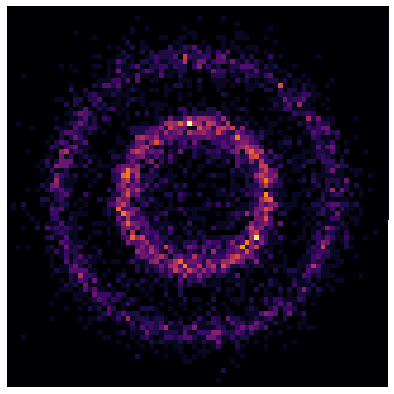}&
        \includegraphics[width=0.17\textwidth]{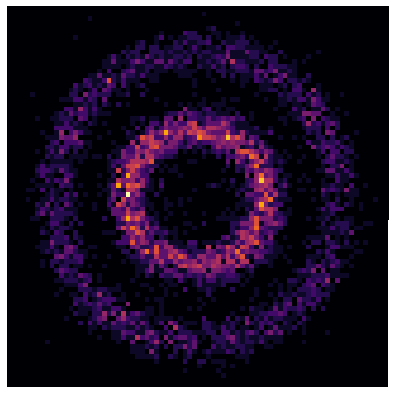}&
        \includegraphics[width=0.17\textwidth]{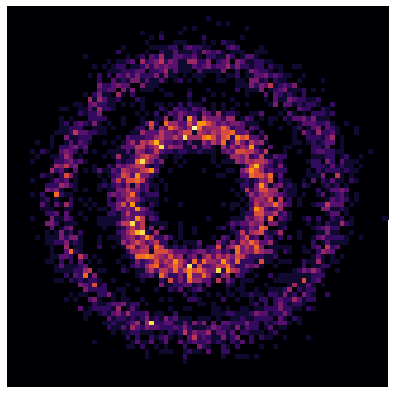}&
        \includegraphics[width=0.17\textwidth]{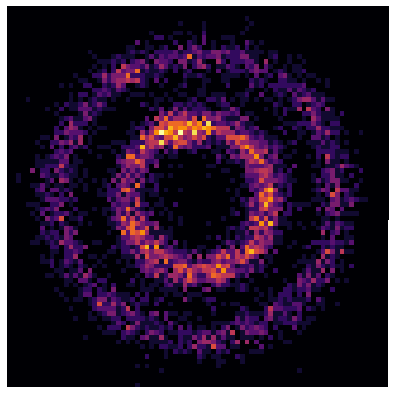}\\
        \\
        \multicolumn{5}{c}{\textbf{MMD$\boldsymbol{^2}$ Values, JKO-Flow (five iterations)}}
        \\
		1.94e-3 & 3.24e-4 & 7.71e-4 & 5.9e-5 & 1.01e-4 \\
		\includegraphics[width=0.17\textwidth]{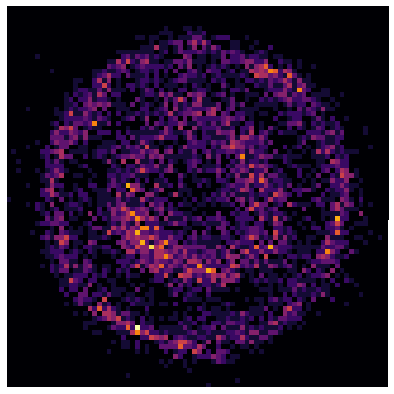}&
        \includegraphics[width=0.17\textwidth]{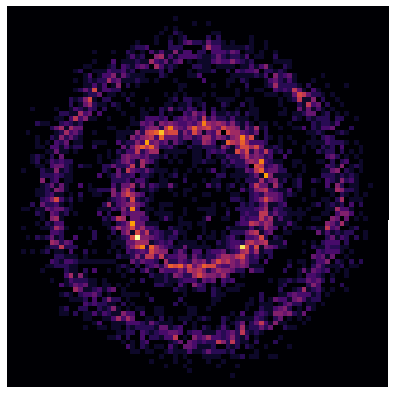}&
        \includegraphics[width=0.17\textwidth]{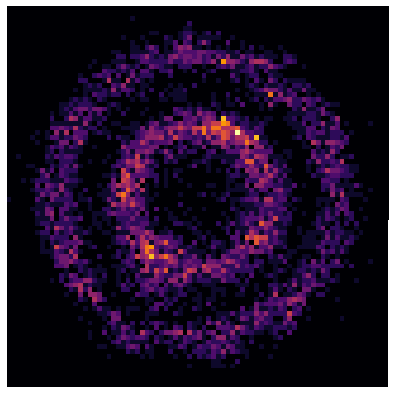}&
        \includegraphics[width=0.17\textwidth]{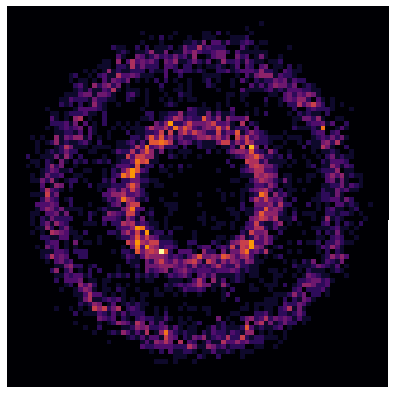}&
        \includegraphics[width=0.17\textwidth]{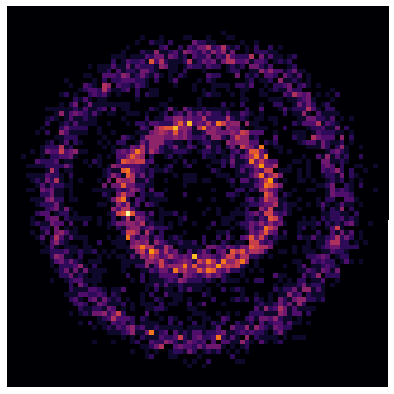} \\
	\end{tabular} 
	\centering 
	\\[2mm]
	\caption{Circles dataset: Generated samples of $\hat{\rho}_0$ using the standard single shot approach (top row). Generated samples using our proposed JKO-Flow using five iterations (bottom row). Here, we fix $\alpha = 5$ and vary the network width $m = 3, 4, 5, 8$, and $16$. JKO-Flow performs competitively even with lower number of parameters.}
	
\end{figure} 

\begin{figure}[ht] 
\centering
	\begin{tabular}{ccccc} 	    
		$m = 3$ & $m = 4$ & $m = 5$  & $m = 8$  & $m = 16$ \\ 
          \\ 
        \multicolumn{5}{c}{\textbf{MMD$\boldsymbol{^2}$ Values, Single Shot}}
		\\
		1.10e-2 & 4.03e-3 & 2.27e-3 & 3.8e-3 & 5.43e-4 \\
		\includegraphics[width=0.17\textwidth]{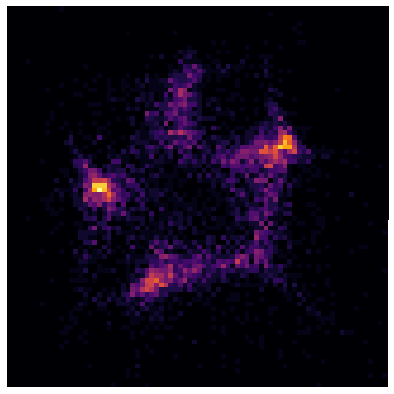}&
        \includegraphics[width=0.17\textwidth]{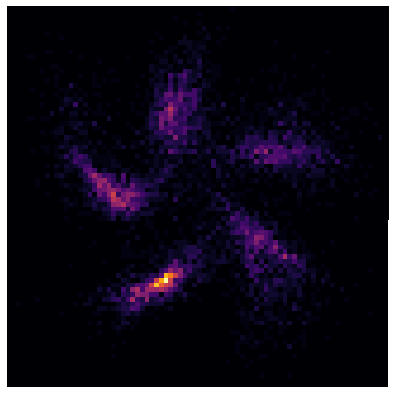}&
        \includegraphics[width=0.17\textwidth]{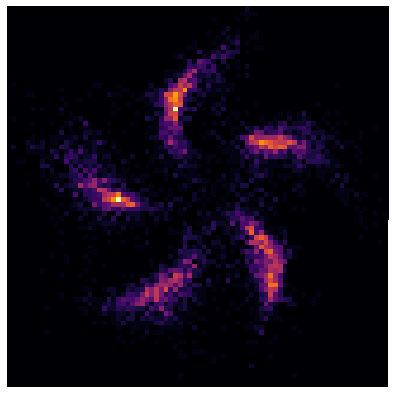}&
        \includegraphics[width=0.17\textwidth]{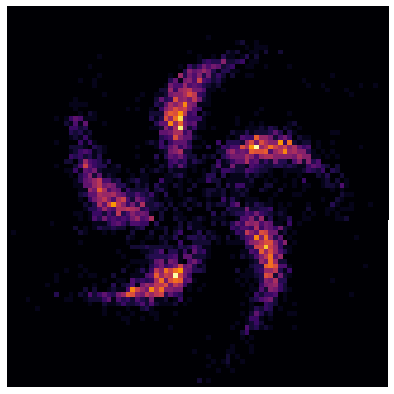}&
        \includegraphics[width=0.17\textwidth]{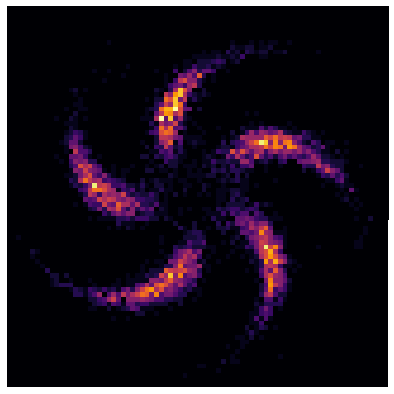}\\
        \\
        \multicolumn{5}{c}{\textbf{MMD$\boldsymbol{^2}$ Values, JKO-Flow (five iterations)}}
        \\
		1.20e-3 & 8.23e-4 & 1.60e-3 & 7e-5 & 2.69e-4 \\
		\includegraphics[width=0.17\textwidth]{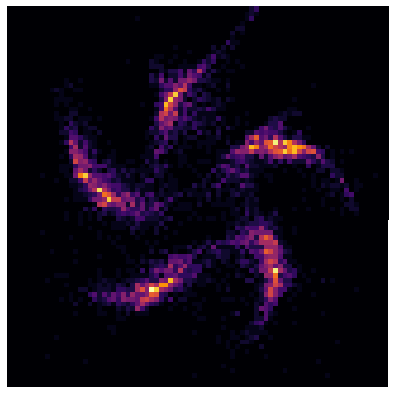}&
        \includegraphics[width=0.17\textwidth]{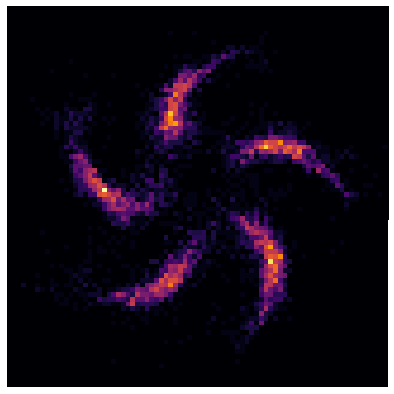}&
        \includegraphics[width=0.17\textwidth]{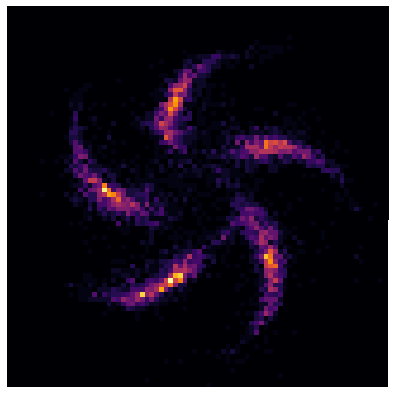}&
        \includegraphics[width=0.17\textwidth]{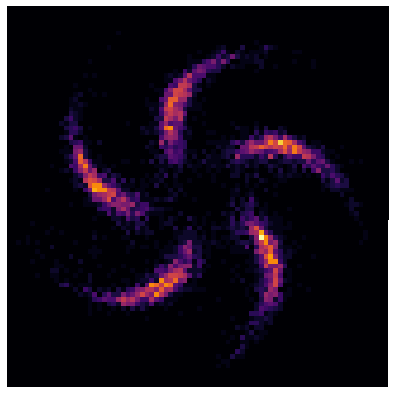}&
        \includegraphics[width=0.17\textwidth]{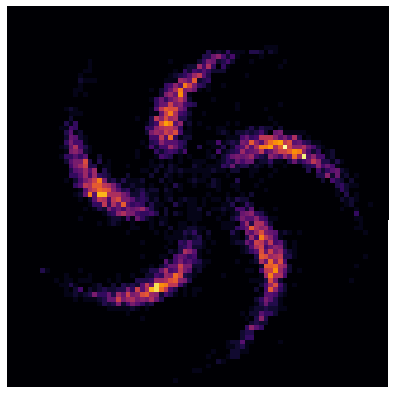} \\
	\end{tabular} 
	\centering 
	\\[2mm]
	\caption{Pinwheel dataset: Generated samples of $\hat{\rho}_0$ using the standard single shot approach (top row). Generated samples using our proposed JKO-Flow using five iterations (bottom row). Here, we fix $\alpha = 5$ and vary the network width $m = 3, 4, 5, 8$, and $16$. JKO-Flow performs competitively even with lower number of parameters.}
	
\end{figure} 
\begin{figure}[ht] 
\centering
	\begin{tabular}{ccccc} 	    
		$m = 3$ & $m = 4$ & $m = 5$  & $m = 8$  & $m = 16$ \\ 
            \\ 
        \multicolumn{5}{c}{\textbf{MMD$\boldsymbol{^2}$ Values, Single Shot}}
		\\
		  5.98e-3 & 4.54e-3 & 5.47e-3 & 1.2e-3 & 3.96e-3 \\
		\includegraphics[width=0.17\textwidth]{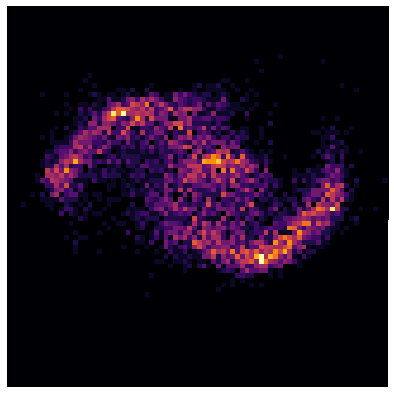}&
        \includegraphics[width=0.17\textwidth]{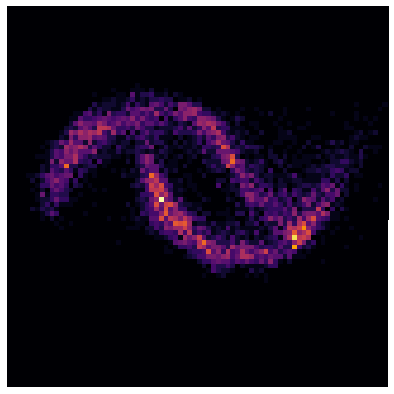}&
        \includegraphics[width=0.17\textwidth]{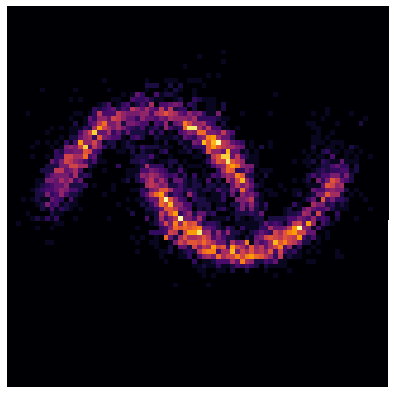}&
        \includegraphics[width=0.17\textwidth]{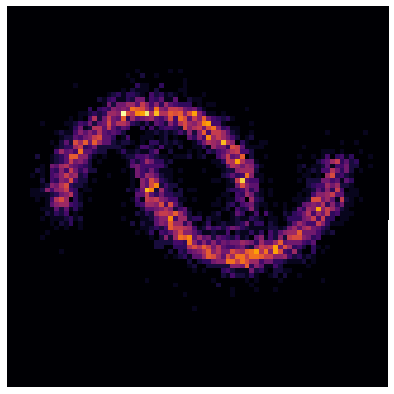}&
        \includegraphics[width=0.17\textwidth]{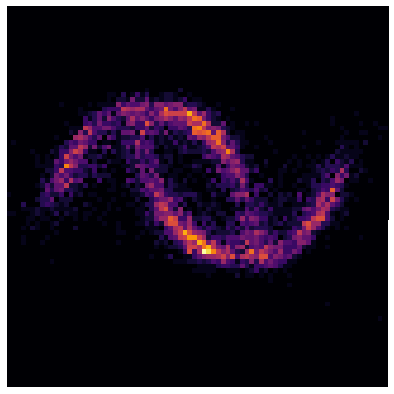}\\
                \\
        \multicolumn{5}{c}{\textbf{MMD$\boldsymbol{^2}$ Values, JKO-Flow (five iterations)}}
        \\
		1.42e-3 & 1.5e-5 & 6.11e-4 & 3.9e-5 & 2.19e-3 \\
		\includegraphics[width=0.17\textwidth]{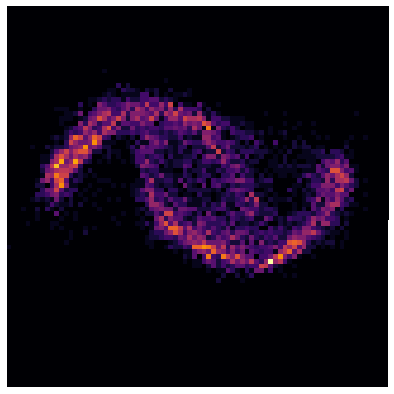}&
        \includegraphics[width=0.17\textwidth]{figures/moons_alpha5_m4_1iter.png}&
        \includegraphics[width=0.17\textwidth]{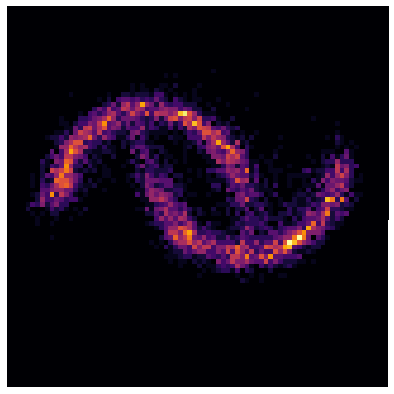}&
        \includegraphics[width=0.17\textwidth]{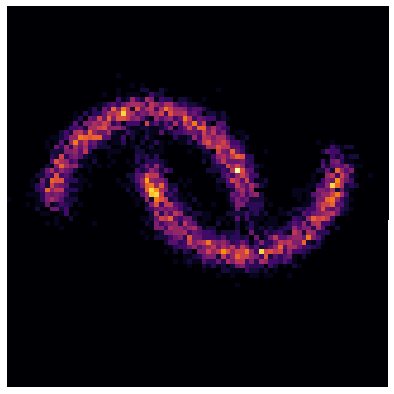}&
        \includegraphics[width=0.17\textwidth]{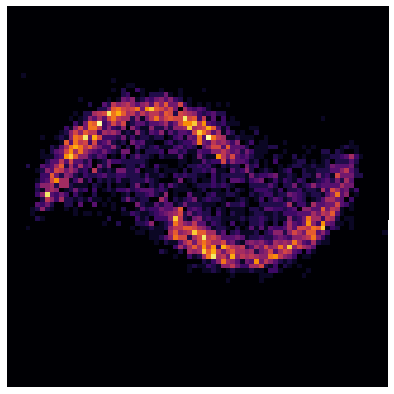} \\
	\end{tabular} 
	\centering 
	\\[2mm]
	\caption{Moons dataset: Generated samples of $\hat{\rho}_0$ using the standard single shot approach (top row). Generated samples using our proposed JKO-Flow using five iterations (bottom row). Here, we fix $\alpha = 5$ and vary the network width $m = 3, 4, 5, 8$, and $16$. JKO-Flow performs competitively even with lower number of parameters.}
\end{figure} 


\begin{figure}[ht]
    \centering
    \addtolength{\tabcolsep}{-6pt} 
    \subfloat[\miniboone{} dimension 16 vs 17\label{fig:miniboone1.1}]{
	   	\begin{tabular}{ccc}
	   	
	   	\textbf{Samples} & Single Shot & \hspace{6pt} JKO Flow \\
	   	$\bfx \sim \rho_0(\bfx)$ & $f(\bfx)$ & \hspace{6pt}$f(\bfx)$ \\
	   	\includegraphics[height=46pt]{figures/data_image.png}
	   	&
	   	\includegraphics[height=46pt]{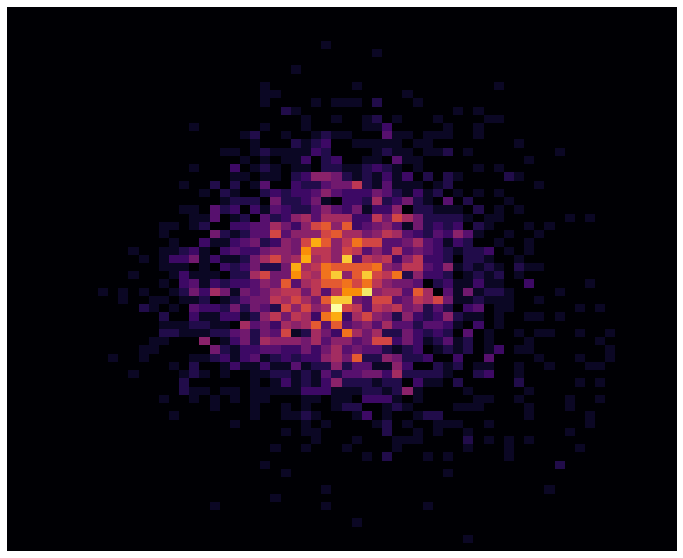}
	   	&
	   	\includegraphics[height=46pt]{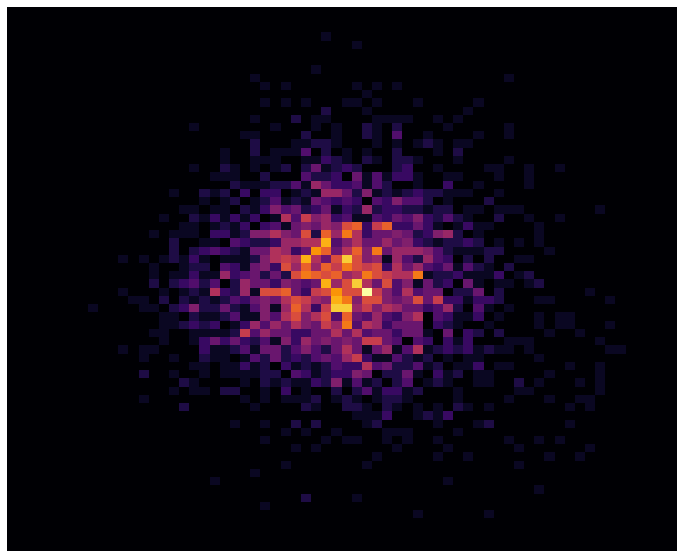} \\
	   	$\bfy \sim \rho_1(\bfy)$ & $f^{-1}(\bfy)$ & \hspace{6pt} $f^{-1}(\bfy)$ \\ 
	   	\includegraphics[height=46pt]{figures/miniboone_y.png}
	   	&
	   	\includegraphics[height=46pt]{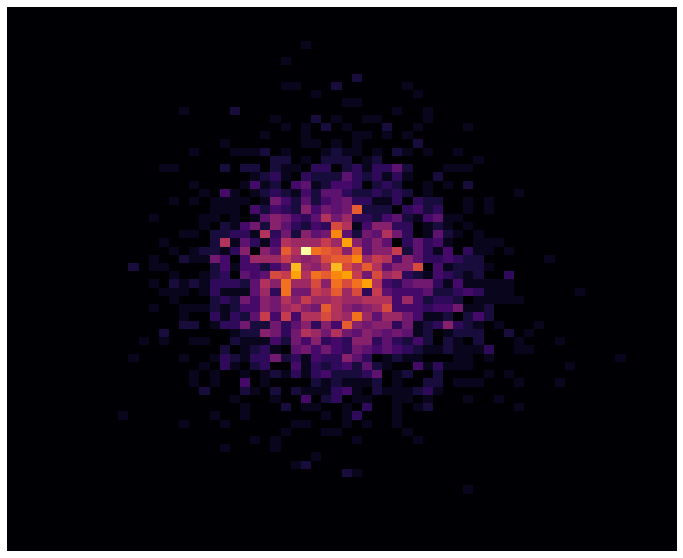}
	   	&
	   	\includegraphics[height=46pt]{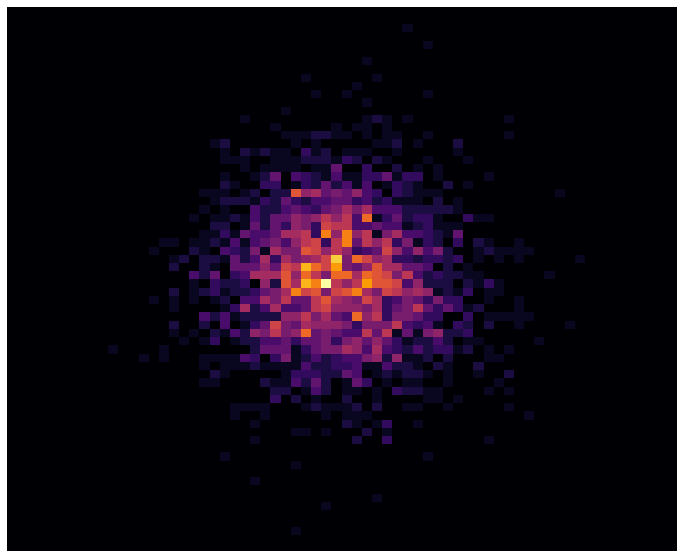}
	   	\end{tabular}%
	}\hspace{10pt}%
  	\subfloat[\miniboone{} dimension 28 vs 29\label{fig:miniboone1.2}]{
	   	\begin{tabular}{ccc}
	   	\textbf{Samples} & Single Shot & \hspace{6pt} JKO Flow \\
	   	$\bfx \sim \rho_0(\bfx)$ & $f(\bfx)$ & \hspace{6pt} $f(\bfx)$ \\
	   	\includegraphics[height=46pt]{figures/data_image_28_v_29.png}
	   	&
	   	\includegraphics[height=46pt]{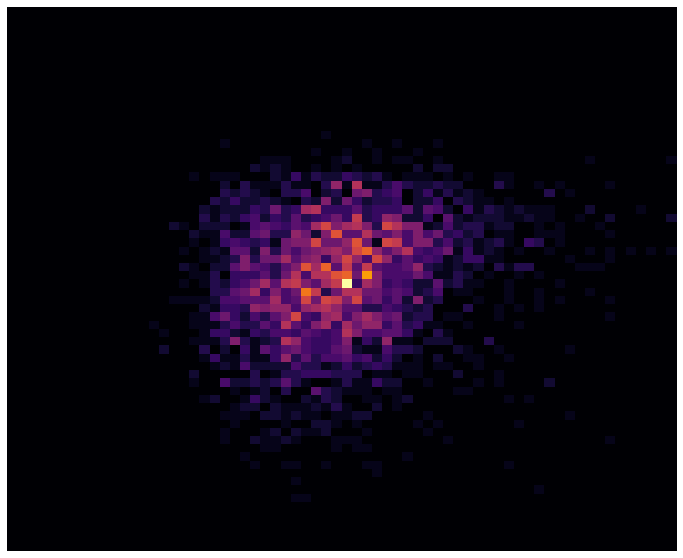}
	   	&
	   	\includegraphics[height=46pt]{figures/miniboone_fx_alpha1_5subproblem_28v29.png} \\
	   	$\bfy \sim \rho_1(\bfy)$ & $f^{-1}(\bfy)$ & \hspace{6pt} $f^{-1}(\bfy)$ \\ 
	   	\includegraphics[height=46pt]{figures/miniboone_y.png}
	   	&
	   	\includegraphics[height=46pt]{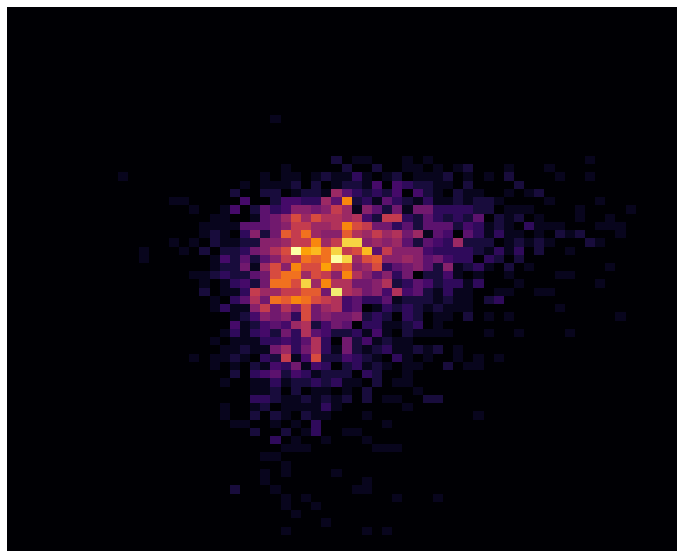}
	   	&
	   	\includegraphics[height=46pt]{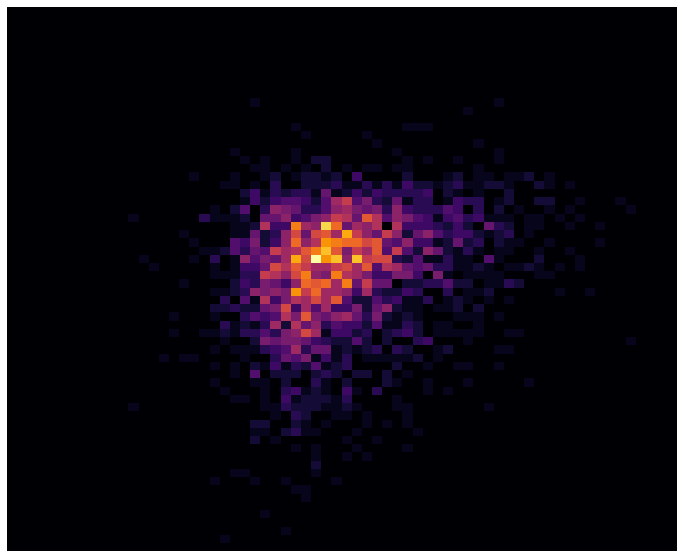}
	   	\end{tabular}%
	}%
    \caption{Generated samples for the 43-dimensional \miniboone{} dataset using the single shot approach and JKO-Flow with 10 iterations for $\alpha = 1$. To visualize the dataset, we show 2-dimensional slices. We show the forward flow $f(x)$ where $x \sim \rho_0$ and the genereated samples $f^{-1}(y)$ where $y \sim \rho_1$.
    }
    \label{fig:miniboone1}
\addtolength{\tabcolsep}{5pt} 
\end{figure}
\begin{figure}[ht]
    \centering
    \addtolength{\tabcolsep}{-6pt} 
    \subfloat[\miniboone{} dimension 16 vs 17\label{fig:miniboone10.1}]{
	   	\begin{tabular}{ccc}
	   	
	   	\textbf{Samples} & Single Shot & \hspace{6pt} JKO Flow \\
	   	$\bfx \sim \rho_0(\bfx)$ & $f(\bfx)$ & \hspace{6pt}$f(\bfx)$ \\
	   	\includegraphics[height=46pt]{figures/data_image.png}
	   	&
	   	\includegraphics[height=46pt]{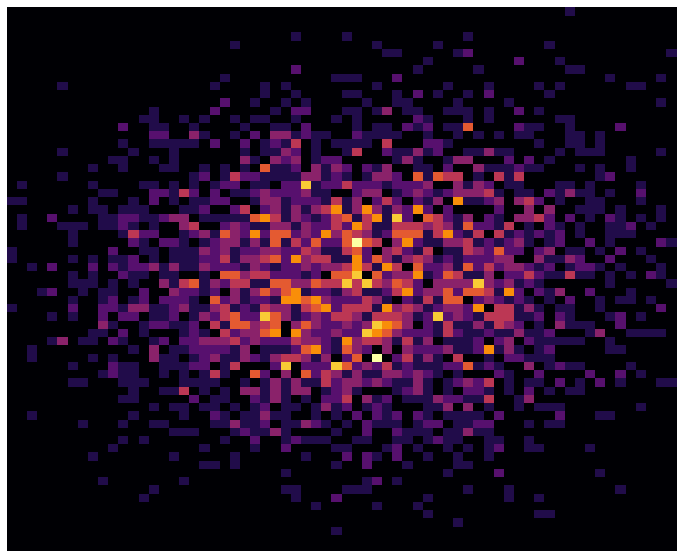}
	   	&
	   	\includegraphics[height=46pt]{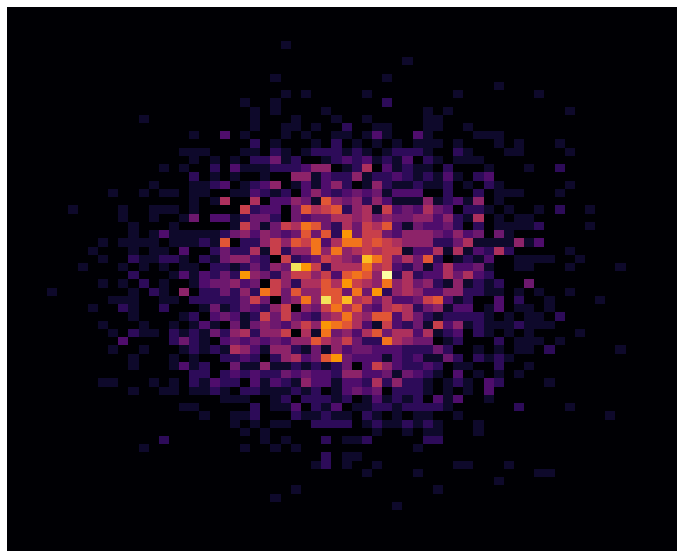} \\
	   	$\bfy \sim \rho_1(\bfy)$ & $f^{-1}(\bfy)$ & \hspace{6pt} $f^{-1}(\bfy)$ \\ 
	   	\includegraphics[height=46pt]{figures/miniboone_y.png}
	   	&
	   	\includegraphics[height=46pt]{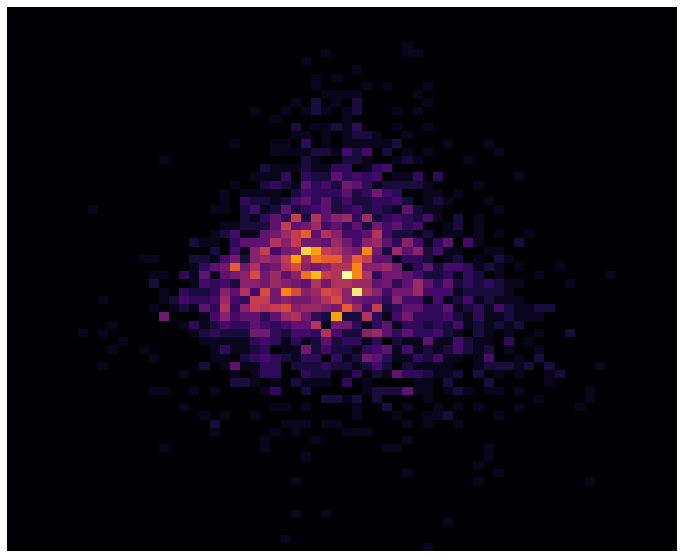}
	   	&
	   	\includegraphics[height=46pt]{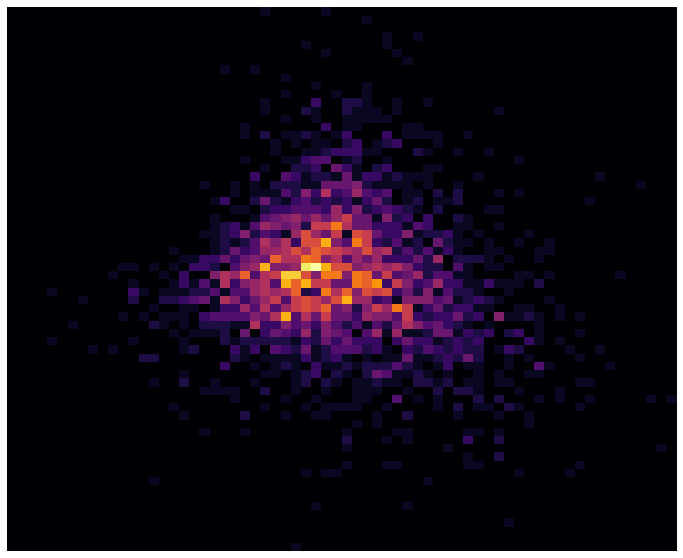}
	   	\end{tabular}%
	}\hspace{10pt}%
  	\subfloat[\miniboone{} dimension 28 vs 29\label{fig:miniboone10.2}]{
	   	\begin{tabular}{ccc}
	   	\textbf{Samples} & Single Shot & \hspace{6pt} JKO Flow \\
	   	$\bfx \sim \rho_0(\bfx)$ & $f(\bfx)$ & \hspace{6pt} $f(\bfx)$ \\
	   	\includegraphics[height=46pt]{figures/data_image_28_v_29.png}
	   	&
	   	\includegraphics[height=46pt]{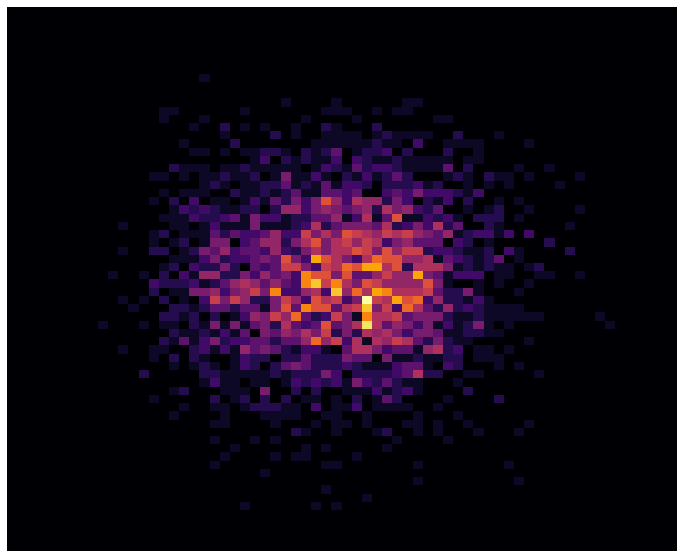}
	   	&
	   	\includegraphics[height=46pt]{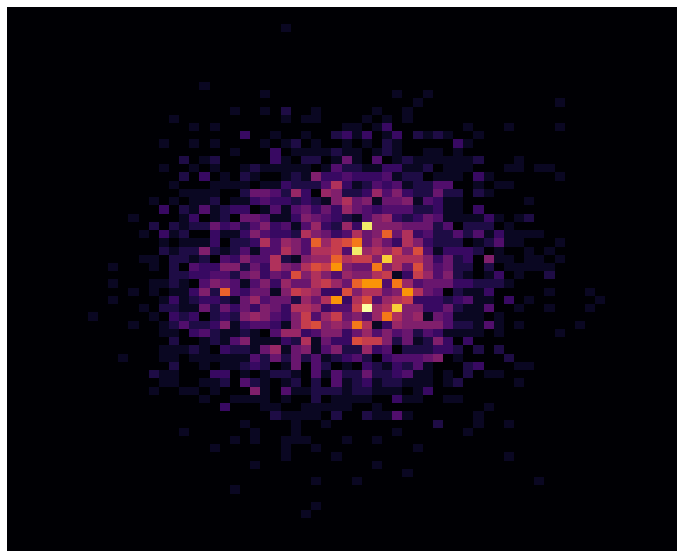} \\
	   	$\bfy \sim \rho_1(\bfy)$ & $f^{-1}(\bfy)$ & \hspace{6pt} $f^{-1}(\bfy)$ \\ 
	   	\includegraphics[height=46pt]{figures/miniboone_y.png}
	   	&
	   	\includegraphics[height=46pt]{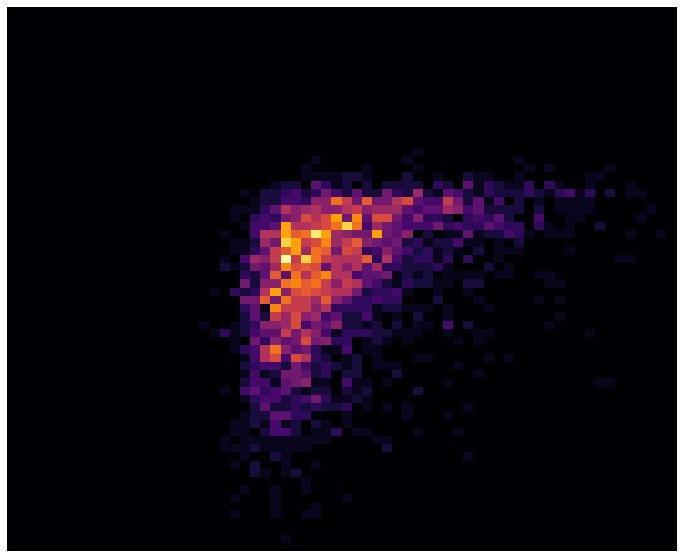}
	   	&
	   	\includegraphics[height=46pt]{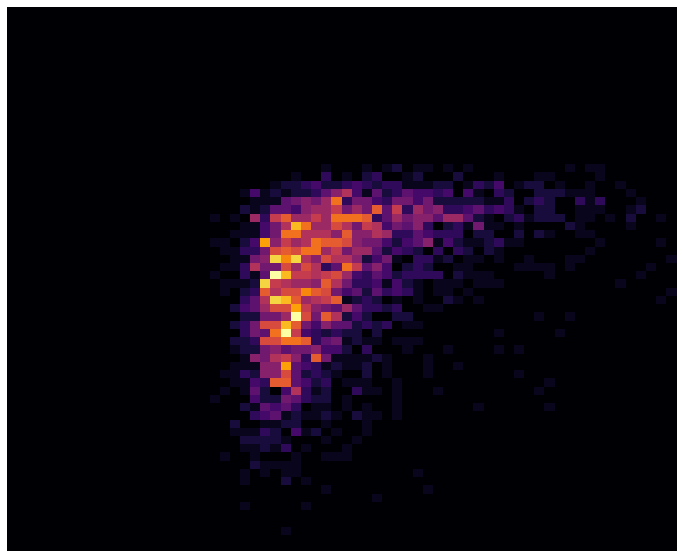}
	   	\end{tabular}%
	}%
    \caption{    Generated samples for the 43-dimensional \miniboone{} dataset using the single shot approach and JKO-Flow with 10 iterations for $\alpha = 10$. To visualize the dataset, we show 2-dimensional slices. We show the forward flow $f(x)$ where $x \sim \rho_0$ and the genereated samples $f^{-1}(y)$ where $y \sim \rho_1$.
    }
    \label{fig:miniboone10}
\addtolength{\tabcolsep}{5pt} 
\end{figure}
\begin{figure}[ht]
    \centering
    \addtolength{\tabcolsep}{-6pt} 
    \subfloat[\miniboone{} dimension 16 vs 17\label{fig:miniboone50.1}]{
	   	\begin{tabular}{ccc}
	   	\textbf{Samples} & Single Shot & \hspace{6pt} JKO Flow \\
	   	$\bfx \sim \rho_0(\bfx)$ & $f(\bfx)$ & \hspace{6pt}$f(\bfx)$ \\
	   	\includegraphics[height=46pt]{figures/data_image.png}
	   	&
	   	\includegraphics[height=46pt]{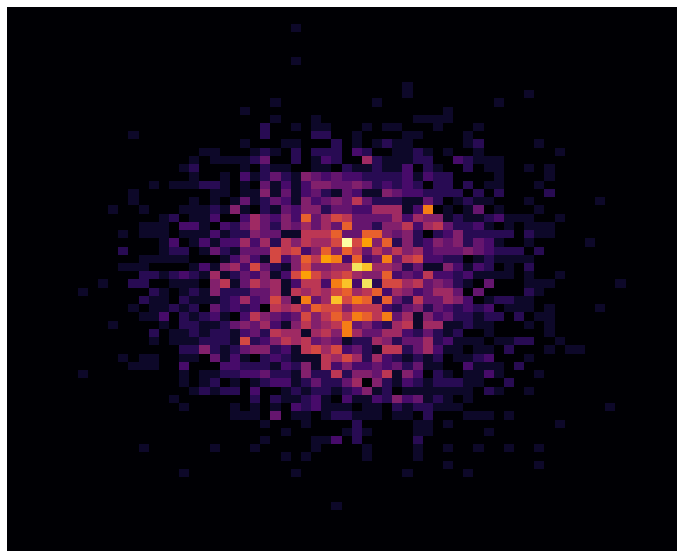}
	   	&
	   	\includegraphics[height=46pt]{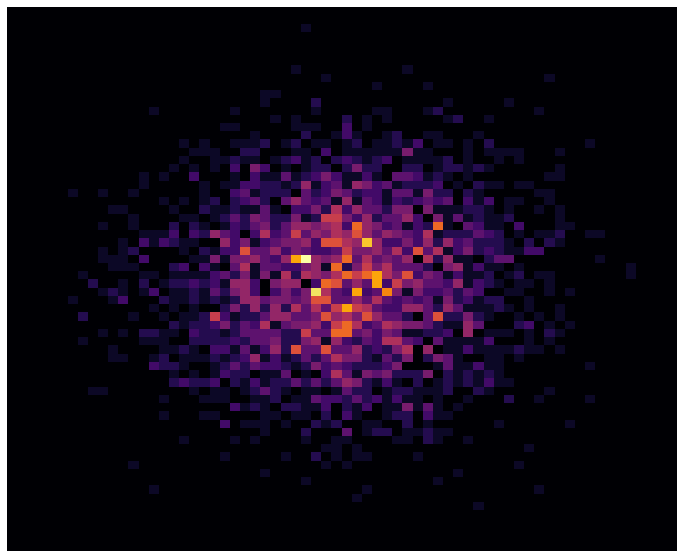} \\
	   	$\bfy \sim \rho_1(\bfy)$ & $f^{-1}(\bfy)$ & \hspace{6pt} $f^{-1}(\bfy)$ \\ 
	   	\includegraphics[height=46pt]{figures/miniboone_y.png}
	   	&
	   	\includegraphics[height=46pt]{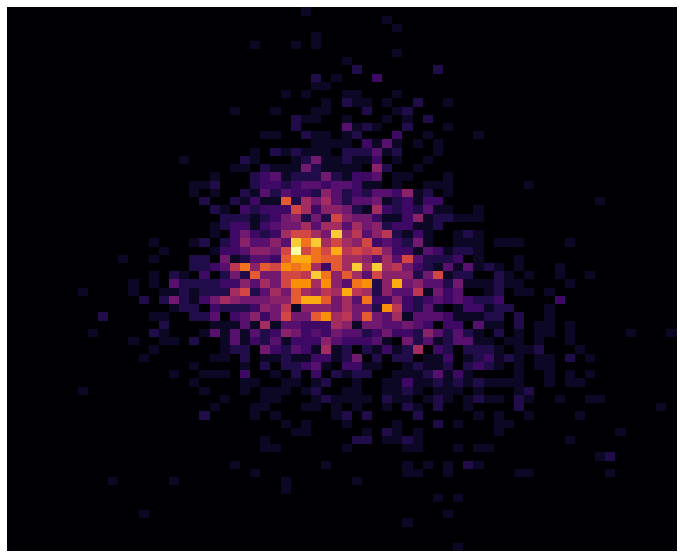}
	   	&
	   	\includegraphics[height=46pt]{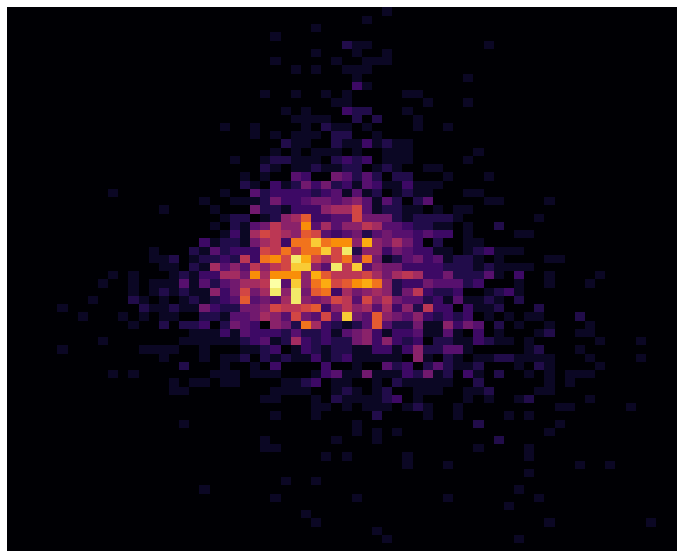}
	   	\end{tabular}%
	}\hspace{10pt}%
  	\subfloat[\miniboone{} dimension 28 vs 29\label{fig:miniboone50.2}]{
	   	\begin{tabular}{ccc}
	   	\textbf{Samples} & Single Shot & \hspace{6pt} JKO Flow \\
	   	$\bfx \sim \rho_0(\bfx)$ & $f(\bfx)$ & \hspace{6pt} $f(\bfx)$ \\
	   	\includegraphics[height=46pt]{figures/data_image_28_v_29.png}
	   	&
	   	\includegraphics[height=46pt]{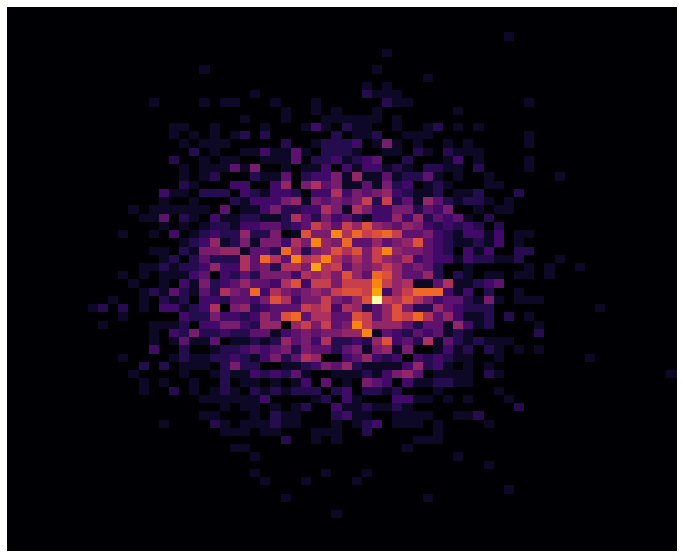}
	   	&
	   	\includegraphics[height=46pt]{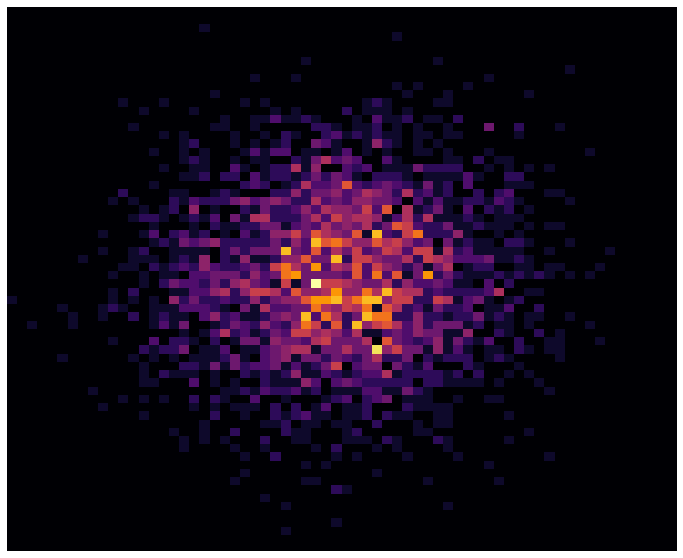} \\
	   	$\bfy \sim \rho_1(\bfy)$ & $f^{-1}(\bfy)$ & \hspace{6pt} $f^{-1}(\bfy)$ \\ 
	   	\includegraphics[height=46pt]{figures/miniboone_y.png}
	   	&
	   	\includegraphics[height=46pt]{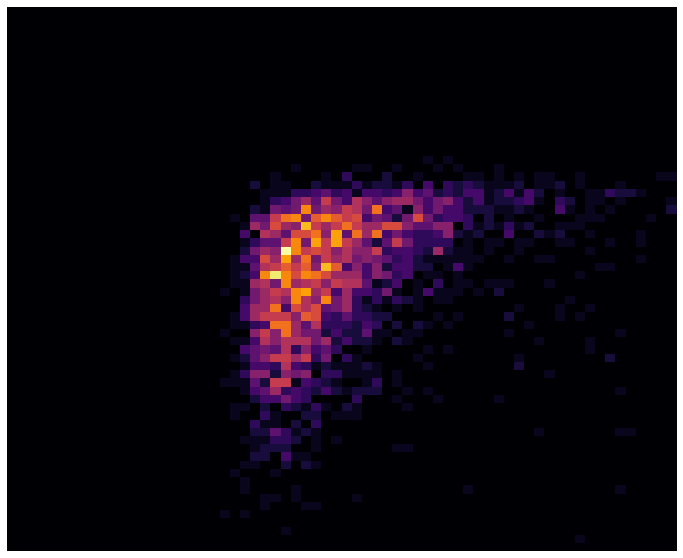}
	   	&
	   	\includegraphics[height=46pt]{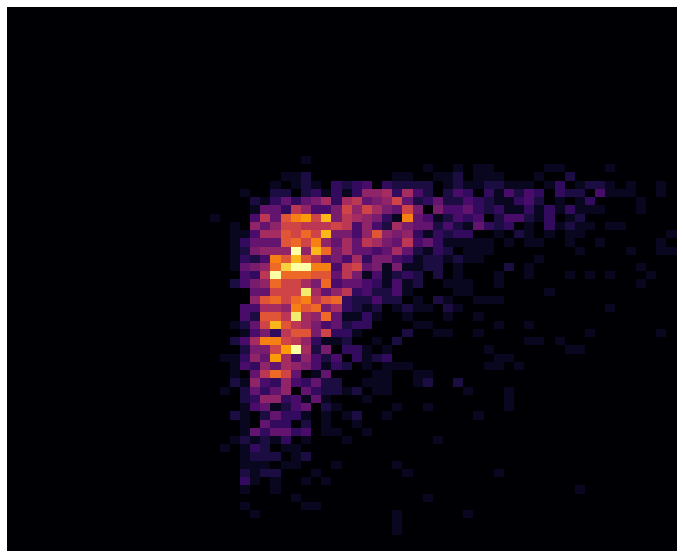}
	   	\end{tabular}%
	}%
    \caption{
    Generated samples for the 43-dimensional \miniboone{} dataset using the single shot approach and JKO-Flow with 10 iterations for $\alpha = 50$. To visualize the dataset, we show 2-dimensional slices. We show the forward flow $f(x)$ where $x \sim \rho_0$ and the genereated samples $f^{-1}(y)$ where $y \sim \rho_1$.
    }
    \label{fig:miniboone50}
\addtolength{\tabcolsep}{5pt} 
\end{figure}

\end{document}